\documentclass[10pt]{article}

 \usepackage{amssymb,amsmath,amsfonts,amssymb}
 \usepackage{graphics,graphicx,color}
 -\marginparwidth \oddsidemargin -4mm \evensidemargin
 \oddsidemargin%
 \usepackage{amsmath}    
 \advance\hoffset by 5mm
 \definecolor{DarkBlue}{rgb}{0.00,0.00,0.55}

 \title{Orientability and energy minimization \\ in liquid crystal models}
 \date{\today}
\author{John M. Ball\footnote{Mathematical Institute, 24-29 St. Giles', Oxford, OX1 3LB, email: ball@maths.ox.ac.uk}  \, and Arghir Zarnescu\footnote{Mathematical Institute, 24-29 St. Giles', Oxford, OX1 3LB, email: zarnescu@maths.ox.ac.uk}}

\newtheorem{definition}{Definition}
\newtheorem{remark}{Remark}
 \newtheorem{lemma} {Lemma}
 \newtheorem{proposition}{Proposition}
 \newtheorem{theorem} {Theorem}

\begin{document}
\maketitle
\begin{abstract}
Uniaxial nematic liquid crystals are  modelled in the Oseen-Frank theory through a unit vector field $n$. This theory has the apparent drawback that it does not respect the head-to-tail symmetry in which $n$ should be equivalent to -$n$. This symmetry is preserved in the constrained Landau~-~de Gennes theory that works with the tensor $Q=s
\big(n\otimes n- \frac{1}{3} Id\big)$.
\par We study the differences and the overlaps between the two theories. These depend on the regularity class used as well as on the topology of the underlying domain. We show that for simply-connected domains and in the natural energy class $W^{1,2}$ the two theories coincide, but otherwise there can be differences between the two theories, which we identify.
\par In the case of planar domains we completely characterise the instances in which the predictions of the constrained Landau~-~de Gennes theory differ from those of the Oseen-Frank theory.
\end{abstract}
\section{Introduction}
\par The challenge of describing  nematic liquid crystals by a model that is both comprehensive and simple enough to manipulate efficiently has led to the existence of several major competing theories. One of the most simple and successful is the  Oseen-Frank theory that describes nematics using  relatively simple tools, namely vector fields, but  has the deficiency of ignoring a physical symmetry of the material. A more complex theory was proposed by de Gennes and  uses matrix-valued functions (Q-tensors). In the simplest constrained case of uniaxial Q-tensors with a constant scalar order parameter these Q-tensors can be interpreted as line fields. We study in this paper the differences between these two theories in this constrained case and establish when the more complicated, but physically more realistic,  theory of de Gennes can be replaced by the simpler Oseen-Frank theory, and when this cannot be done. Of particular interest to our study are  liquid crystal `defects' and we examine the instances when one would detect fake `defects' in the material resulting from using the more simplistic Oseen-Frank model instead of the constrained Landau~-~de Gennes model (we refer to such a situation as having a {\it non-orientable line field}).

In nematic liquid crystals  the rod-like
molecules tend to align, locally,  along a preferred direction. This
is modeled by assigning at each material point $x$ in the region $\Omega$ occupied
by the liquid crystal   a probability measure  $\mu(x,\cdot):\mathcal{L}(\mathbb{S}^2)\to\mathbb
[0,1]$ for describing the orientation of the molecules, where $\mathcal{L}(\mathbb{S}^2)$ denotes the family of
Lebesgue measurable sets on the unit sphere. Thus $\mu(x,A)$ gives
the probability that the molecules with centre of mass in a very small neighbourhood of  the point
$x\in\Omega$ are pointing in a direction contained in $A\subset
\mathbb{S}^2$.
\par Nematic liquid crystals are locally
invariant with respect to reflection in the plane perpendicular to
the preferred direction. This is commonly referred to as the
`head-to-tail' symmetry, see \cite{gennes}. This means that
$\mu(x,A)=\mu(x,-A),\mbox{ for all }
x\in\Omega,A\subset\mathcal{L}(\mathbb{S}^2)$. Note that because of
this symmetry the first moment of the probability measure vanishes:
$$\langle
p\rangle\stackrel{\rm{def}}{=}\int_{\mathbb{S}^2}p\,d\mu(p)=\frac{1}{2}\left[\int_{\mathbb{S}^2}p\,d\mu(p)+
\int_{\mathbb{S}^2}-p\,d\mu(-p)\right]=0.$$
\par The first nontrivial information on $\mu$ comes from the
tensor of second moments:
$$M_{ij}\stackrel{\rm{def}}{=}\int_{\mathbb{S}^2}p_ip_j\,d\mu(p),\,i,j=1,2,3.$$
We have $M=M^T$ and
$\textrm{tr}\,M=\int_{\mathbb{S}^2}d\mu(p)=1$. Let $e$ be a unit
vector. Then
$$e\cdot Me=\int_{\mathbb{S}^2}(e\cdot p)^2\,d\mu(p)=\langle
\cos^2(\theta)\rangle$$ where $\theta$ is the angle between $p$ and
$e$.
\par If the orientation of the molecules is equally distributed in
all directions we say that the distribution is {\it isotropic} and
then $\mu=\mu_0$ where $d\mu_0(p)=\frac{1}{4\pi}dA$. The
corresponding second moment tensor is
$$M_0\stackrel{\rm{def}}{=}\frac{1}{4\pi}\int_{\mathbb{S}^2} p\otimes
p\,dA=\frac{1}{3}Id$$ (since $\int_{\mathbb{S}^2}
p_1p_2\,d\mu(p)=0,\,\int_{\mathbb{S}^2}p_1^2\,d\mu(p)=\int_{\mathbb{S}^2}p_2^2\,d\mu(p)=\int_{\mathbb{S}^2}p_3^2\,d\mu(p)$
and $\textrm{tr}\,M_0=1$).
\par The de Gennes order-parameter tensor $Q$ is defined as

\begin{equation}
Q\stackrel{\rm{def}}{=}M-M_0=\int_{\mathbb{S}^2}\left(p\otimes
p-\frac{1}{3}Id\right)\,d\mu(p)\label{defq}
\end{equation} and measures the deviation of the second moment
tensor from its isotropic value.
\par Since $Q$ is symmetric and $\textrm{tr}\,Q=0$, $Q$ has, by the
spectral theorem, the representation:

$$Q=\lambda_1\hat e_1\otimes\hat e_1+\lambda_2\hat e_2\otimes\hat
e_2-(\lambda_1+\lambda_2)\hat e_3\otimes \hat e_3$$ where $\hat
e_1,\,\hat e_2,\,\hat e_3$ is an orthonormal basis of eigenvectors
of $Q$ with the corresponding eigenvalues $\lambda_1,\lambda_2$ and
$\lambda_3=-(\lambda_1+\lambda_2)$.
\par When two of the eigenvalues are equal (and non-zero) the order parameter $Q$ is
called {\it uniaxial}, otherwise being {\it biaxial}. The case $Q=0$ is the {\it isotropic} state. Equilibrium configurations of liquid crystals are obtained, for instance, as
energy minimizers, subject to suitable boundary conditions. The simplest commonly used energy functional is
\begin{equation}
\label{energy}
 \mathcal{F}_{LG}[Q]=\int_{\Omega}\left[ \frac{L}{2}\sum_{i,j,k=1}^3 Q_{ij,k}Q_{ij,k}-\frac{a}{2}\textrm{tr}\,Q^2 -
\frac{ b}{3}\textrm{tr}\,Q^3 +\frac{ c}{4}\left(\textrm{tr}\,Q^2\right)^2\right]\,dx
\end{equation} where $a,b,c$ are temperature and material dependent constants and $L>0$ is the elastic constant. In the physically significant limit $L\to 0$ (and for appropriate boundary conditions) we have that the energy minimizers are suitably approximated by minimizers of  the corresponding {\it `Oseen-Frank energy functional'}
\begin{displaymath}
\mathcal{F}_{OF}[Q]=\int_{\Omega} \sum_{i,j,k=1}^3 Q_{ij,k}Q_{ij,k}\,dx
\end{displaymath} in the restricted class of $Q\in W^{1,2}$, with $Q$ uniaxial a.e., so that

\begin{equation}
Q=s\left(n\otimes n-\frac{1}{3}Id\right) \label{q}
\end{equation} with $s\in\mathbb{R}$ (an explicit constant depending on $a,b$ and $c$) and $n(x)\in\mathbb{S}^2$ a.e. $x\in\Omega$.

\par This convergence, as $L\to 0$, was studied initially in \cite{mz} and further refined in \cite{nz}. In the following we will restrict ourselves to studying tensors $Q$ that admit a representation as in (\ref{q}) and we will further refer to this as {\it the constrained Landau~-~de Gennes theory}.
\par Taking into account the definition  (\ref{defq}) of the tensor $Q$ we have that
\begin{displaymath}
Qn\cdot n=\frac{2}{3}s=\langle (p\cdot
n)^2-\frac{1}{3}\rangle=\langle \cos^2\theta-\frac{1}{3}\rangle
\end{displaymath} where $\theta$ is the angle between $p$ and $n$. Hence
$s=\frac{3}{2}\langle \cos^2\theta-\frac{1}{3}\rangle$ and so we
must necessarily have $-\frac{1}{2}\le s\le 1$ with $s=1$ when we
have perfect ordering parallel to $n$, $s=-\frac{1}{2}$ when all
molecules are perpendicular to $n$ and $s=0$ iff  $Q=0$ (which
occurs if $\mu$ is isotropic). Thus $s$ is a measure of order and is
called the {\it scalar order parameter associated to the tensor}
$Q$. In the physical literature it is often assumed that $s$ is
constant almost everywhere. For experimental determinations of $s$
see for instance \cite{crawford}.

\par We continue working under this assumption, that $s$ is a
non-zero constant, independent of $x\in\Omega$. Thus, taking into
account the representation (\ref{q}) of the Q-tensor  we have that,
for constant $s$, there exists a bijective correspondence between
$Q$ that have the representation (\ref{q}) and pairs $\{n,-n\}$
with $n\in\mathbb{S}^2$. Hence we can think of $Q$ as the line
joining $n$ and $-n$. We can thus identify the space of $Q$ as in
(\ref{q}) with the real projective space $\mathbb{R}P^2$; see the
beginning of Section $3$ for more details.

\par Traditionally in the mathematical modeling of liquid crystals the
 Q-tensor (\ref{q}) has been replaced by an  {\it oriented} line, a line
with a direction, i.e. $ n\in \mathbb{S}^2$. This is done in the
Oseen-Frank theory, \cite{oseenfrank}.
\par In this paper we analyze the consequences of this assumption,  from the
point of view of energy minimization. We show that taking into
account the {\it possible} unorientedness of the locally preferred
direction produces a theory that includes the traditional approach
(an oriented vector field has an unoriented counterpart but not
necessarily the other way around) and exhibits  features not seen by
the traditional approaches.
\par We analyze first the relation between line  and vector fields by
determining when a line field can be `oriented' into a vector
field, {\it without changing its regularity class}. We show that this
depends both on the regularity of the line field and the topology of
the domain. We also show that, under suitable assumptions, the
orientability of a line field on a $2D$ bounded domain can be determined
just by knowing the orientability of the line field on the boundary
of the domain.
\par When the class of line fields is strictly larger then the class
of vector fields it is possible that a global energy minimizer is
a line field that is `non-orientable' i.e. cannot be reduced to a
vector field. These are precisely the instances in which the
traditional, Oseen-Frank, theory would fail to recognize a global
minimizer, and would indicate as a global minimizer one that is
physically just a local minimizer, a minimizer in the class of
orientable line fields.  We analyze this situation in the case of
planar line fields  for a simple energy functional, and provide
necessary and sufficient conditions for the global energy minimizer
to be non-orientable, in terms of an integer programming problem.
\par The paper is organized as follows: in Section $2$ we
define rigorously the notions  used, in particular orientability in
the Sobolev space setting, and study the relation between the
regularity of a line field and that of the corresponding vector
field (when it exists). In Section  $3$ we determine the
conditions under which a line field can be oriented into a vector
field and show that on $2D$ domains orientability can be checked at the boundary. In
 Section $4$ we provide analytic orientability criteria on $2D$ domains; that
is we show how to reduce the topological problem of checking the
orientability of the line fields to an analytic computation. In Section $5$ we study the
orientability of the global energy minimizers for planar line
fields.

\bigskip
\section{Preliminaries and notation}
\label{preliminaries}

\par For the rest of the paper  we fix $s\in [-\frac{1}{2},1]$ to be a given non-zero constant and
 $Q$ will denote a $3\times 3$ traceless and symmetric matrix that admits the representation:
\begin{equation}
Q=s\left(n\otimes n-\frac{1}{3}Id\right) \label{Q}
\end{equation} where $n\in\mathbb{R}^3, |n|=1$ and $s$ is the given
non-zero constant.

\par We denote by $\mathcal{M}^{3\times 3}(\mathbb{R})$ the set of $3\times 3$ matrices with real values. In general for an arbitrary symmetric and traceless  $Q\in\mathcal{M}^{3\times 3}(\mathbb{R})$ there might  be no
$n\in\mathbb{S}^2$ so that $Q$ has a representation as in (\ref{Q}).
The necessary and sufficient conditions for a $3\times 3$ matrix to
have such a representation are:

\begin{proposition} For fixed $s\in\mathbb{R}$  and a matrix $Q\in\mathcal{M}^{3\times
3}(\mathbb{R})$ that is symmetric with trace zero the following conditions are equivalent:

\par $($i$)$ $Q=s(n\otimes n-\frac{1}{3}Id)$ for some $n\in\mathbb{S}^2$
\par $($ii$)$ $Q$ has two equal eigenvalues equal to $-\frac{s}{3}$.
\par $($iii$)$ $\det\,Q=-\frac{2s^3}{27}\, ,\,\rm{tr}\,Q^2=\frac{2s^2}{3}.$
\end{proposition}
\smallskip
\begin{remark} From $(iii)$ it follows that a necessary and sufficient condition for $Q$ to have the represenation $(i)$ for some $s\in\mathbb{R}$ is that $\left(\textrm{tr}\,Q^2\right)^3=54\,(\det\,Q)^2$.
\end{remark}
\smallskip
{\bf Proof.} By the spectral theorem $(ii)$ holds if and only if
\begin{displaymath}
Q=-\frac{s}{3}\left(Id-n\otimes n\right)+\frac{2s}{3}n\otimes n
\end{displaymath} for some $n\in\mathbb{S}^2$, which is $(i)$.

\par If the eigenvalues of $Q$ are $\lambda_1,\lambda_2,\lambda_3$, so that $\textrm{tr}\,Q=\lambda_1+\lambda_2+\lambda_3=0$, the characteristic equation for $Q$ is

\begin{displaymath}
\det(\lambda Id-Q)=(\lambda-\lambda_1)(\lambda-\lambda_2)(\lambda-\lambda_3)=\lambda^3-\frac{1}{2}\textrm{tr}(Q^2)\lambda-\det Q=0
\end{displaymath} since $\lambda_1\lambda_2+\lambda_2\lambda_3+\lambda_3\lambda_1=\frac{1}{2}[(\lambda_1+\lambda_2+\lambda_3)^2-(\lambda_1^2+\lambda_2^2+\lambda_3^2)]$. So the eigenvalues are $-\frac{s}{3},-\frac{s}{3},\frac{2s}{3}$ if and only if $(iii)$ holds. $\Box$

\smallskip\par For the fixed non-zero constant $s\in
[-\frac{1}{2},1]$  we denote

\begin{equation}
\boldsymbol{\mathcal{Q}}=\left\{Q\in \mathcal{M}^{3\times
3}(\mathbb{R});Q=s\left(n\otimes n-\frac{1}{3}Id\right)\,\textrm{
for some } n\in\mathbb{S}^2\right\}.
\end{equation}
Also, let us denote

\begin{equation}
\boldsymbol{\mathcal{Q}}_2=\left\{Q\in \boldsymbol{\mathcal{Q}};
\,\,Q=s\left((n_1,n_2,0)\otimes
(n_1,n_2,0)-\frac{1}{3}Id\right)\right\}
\label{q2}
\end{equation} corresponding
to the `planar' unit vectors $n=(n_1,n_2,0)$.

\par Note that the bijective identification of $\mathbb{R}P^2$ with
$\boldsymbol{\mathcal{Q}}$ allows one to endow
$\boldsymbol{\mathcal{Q}}$ with a structure of
 a two-dimensional manifold. Similarly $\boldsymbol{\mathcal{Q}}_2$ can be
bijectively identified with $\mathbb{R}P^1$, and thus can be given a
structure of a  one-dimensional manifold; see the beginning of
Section $3$ for details.

\par We define the projection operator
$P:\mathbb{S}^2\to \boldsymbol{\mathcal{Q}}$ by:
\begin{equation}
P(n)\stackrel{\rm{def}}{=}s\left(n\otimes n-\frac{1}{3}Id\right). \label{p}
\end{equation}
 Note that one has $P(n)=P(-n)$.  Thus the operator $P$
 provides us with a way of `unorienting' an $\mathbb{S}^2$-valued vector field.
\begin{figure}
 \centering
\includegraphics[scale=0.4]{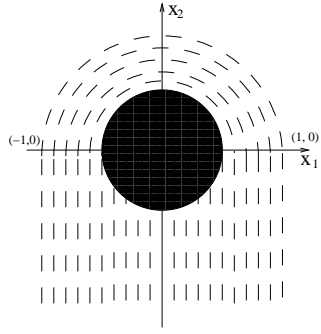}
\caption{A non-orientable director field} \label{canonicalexample}
\end{figure}

 \par  In order  to know when the study of unoriented director fields can be reduced
 to that of oriented ones one needs to know when the opposite is possible, that is to
`orient' a $\boldsymbol{\mathcal{Q}}$-valued vector field. This
should be done without creating `artificial defects', that is discontinuities
of the vector field that were not present in the line field.
Orienting a line field  is sometimes impossible as can be seen by
the example in the Figure $1$ where a line field is defined outside
of the circle centered at $(0,0)$ of radius $\frac{1}{2}$.
Heuristically one sees that by trying to orient the line field for
$x_2>0$ and then for $x_2<0$ one would get a discontinuity. A
rigorous proof of the non-orientability of the line field in Figure
$1$ will be provided in Section $4$, Lemma ~\ref{lemma:specEX}.

\par We define an open set in $\mathbb{R}^m$ to be of class $C^k, k\ge 0$, respectively Lipschitz,  if for any
point $x\in\partial\Omega$ there exist a $\delta>0$ and an
orthonormal coordinate system $Y=(y',y_m)=(y_1,y_2,\dots,y_m)$  with origin at $x$
together with a function of class $C^k$ (respectively Lipschitz)
$f:\mathbb{R}^{m-1}\to\mathbb{R}$ such that
\begin{equation}
\mathcal{U}_{\delta}:=\{y\in\Omega:|y|<\delta\}=\{y\in\mathbb{R}^m:\,y_m>f(y'),\,|y|<\delta\}.
\label{u}
\end{equation}
\par This definition allows one  to consider the open sets of class $C^k$ (respectively Lipschitz) as manifolds with boundary and ensures that the topological boundary of the set coincides with the boundary of the manifold with boundary (see the discussion in the Appendix for more details).
\par In the rest of the paper  an open and connected set is called a {\it domain}.
\par Let us define, in a standard manner (see also
\cite{HALI03} and the references there)  the Sobolev spaces  $W^{1,p}(M,N)$ where $M$ is a $C^0$ manifold and $N$ is a smooth manifold isometrically
embedded in $\mathbb{R}^l$:
$$W^{1,p}(M,N)=\{f:f\in W^{1,p}(M,\mathbb{R}^l), f(x)\in N,\,a.e.\, x\in M\}.$$
In the case when $M$ is a Lipschitz manifold with non-empty boundary, we define
(see also \cite{HALI05}):
$$W^{1,p}_{\varphi}(M,N)=\{f:f\in W^{1,p}(M,\mathbb{R}^l), f(x)\in N,\,a.e.\, x\in M,\,{\rm{Tr}}\,f =\varphi\}$$
where $\rm{Tr}$ denotes the trace operator (\cite{ADAMS75}).
\par We will also sometimes denote the Sobolev space $W^{1,p}(M,N)$, respectively the $L^p$ space $L^p(M,N)$, just by
$W^{1,p}(M)$, respectively $L^p(M)$, when it is clear from the context what the target space is.

\par We define now  orientability for line fields in a Sobolev
space:
\begin{definition}
Let $\Omega\subset \mathbb{R}^d$ be a domain. We say that $Q\in
W^{1,p}(\Omega,\boldsymbol{\mathcal{Q}})$ is  {\it orientable} if there
exists an $n\in W^{1,p}(\Omega,\mathbb{S}^2)$ such that $p(n)=Q$
$\mathcal{L}^d$ a.e. Otherwise we call $Q$   {\it non-orientable}.
\end{definition}

\par We show first that if a line field is orientable there can be
just two orientations, that differ  by change of sign.

\begin{proposition}  Let $\Omega\subset\mathbb{R}^3$ be a domain. An orientable line field $Q(x)\in W^{1,p}(\Omega,\boldsymbol{\mathcal{Q}}),\,1\le p\le \infty$, can
 have only two orientations. More precisely, if $n,m\in
 W^{1,p}(\Omega,\mathbb{S}^2)$ with $P(n)=P(m)=Q$  and $n\not=m$ a.e. we have that
 $m=-n$ a.e. in $\Omega$.
 \label{prop:2orient}
 \end{proposition}

 \smallskip\noindent {\bf Proof.} We have that $m=\tau n\in
 W^{1,p}(\Omega,\mathbb{S}^2)$ with $\tau^2(x)=1$ a.e. In order to
 obtain the conclusion it suffices to show that $\tau(x)$ has
 constant sign almost everywhere.
  \par For a.e. $x_2,x_3$ we have that $\tau(x)n(x)$  has a representative that is absolutely
 continuous in $x_1$ and also $n(x)$ has a representative that is absolutely continuous in
 $x_1$. Hence $\tau(x)n(x)\cdot n(x)=\tau(x)$ is  absolutely continuous in
 $x_1$. A similar statement holds for $x_2$ and $x_3$.

 \par For any $\bar x=(\bar x_1,\bar x_2,\bar x_3)\in \Omega$ let  $\bar\varepsilon>0$  be some number so that $B_{\bar\varepsilon}(\bar x)\subset\Omega$.  On the intersection of almost any line $x_1=const$ with the ball $B_{\bar\varepsilon}(\bar x)$ we have that $\tau(x)$ is  absolutely continuous and, since $\tau^2(x)=1$, we get that $\tau$ is constant.  We take an arbitrary $\varphi\in C_0^\infty(B_{\bar\varepsilon}(\bar x))$ and using Fubini's theorem we obtain:

 \begin{equation}
 \int_{B_{\bar \varepsilon}(\bar x)}\tau(x)\varphi_{,1}(x)\,dx=\int_{B_{\bar\varepsilon}(\bar x)}(\tau\varphi)_{,1}(x)\,dx=0.
 \end{equation}

 \par We obtain thus that the weak derivative  $\tau_{,1}$ is zero in $B_{\bar\varepsilon}(\bar x)$. Similarily
 $\tau_{,2},\tau_{,3}$  are zero  in $B_{\bar\varepsilon}(\bar x)$. Thus $\nabla \tau=0$
 and $\tau$ is constant in $B_{\bar\varepsilon}(\bar x)$. The connectedness of $\Omega$ implies that $\tau$ is constant in $\Omega$. Hence $\tau\equiv 1$ a.e. in $\Omega$ or
 $\tau\equiv-1$ a.e. in $\Omega$. $\Box$

\par Orientability in the open set $\Omega$ implies orientability at the
boundary:
\begin{proposition}
Let $\Omega\subset\mathbb{R}^3$ be a bounded open set with Lipschitz
boundary. Let $Q\in W^{1,p}(\Omega,\boldsymbol{\mathcal{Q}})$ be
orientable, so that $Q=s\left(n\otimes n-\frac{1}{3}Id\right)$ a.e. in $\Omega$ for
some $n\in W^{1,p}(\Omega,\mathbb{S}^2)$. Then $\rm{Tr}\,n\in L^p(\partial\Omega,\mathbb{S}^2)$ and  $Q$ is orientable on
the boundary, i.e.

\begin{equation}
\rm{Tr}\,Q=s\left(\rm{Tr}\,n\otimes\rm{ Tr}\,n-\frac{1}{3}Id\right)
\label{eq:boundaryorient}
\end{equation} with $\rm{Tr}\,Q\in L^p(\partial\Omega,\boldsymbol{\mathcal{Q}})$.
\label{boundaryorient}
\end{proposition}

\smallskip\noindent{\bf Proof.} Choose some open ball $B$ containing
$\bar\Omega$. We can extend $n$ to an $\tilde n\in
W^{1,p}(B,\mathbb{R}^3)$. Truncate each component $\tilde n_i$ by
$1$, i.e. define
\begin{equation}
\bar n_i(x)=\left\{\begin{array}{ll}
            1  & \textrm{ if $\tilde n_i(x)\ge 1$}\\
            \tilde n_i(x) &\textrm{ if  $|\tilde n_i(x)|<1$}\\
            -1 &\textrm{ if  $\tilde n_i(x)\le-1$}
            \end{array}\right. .
\end{equation}
\par Then $\bar n\in L^\infty(B,\mathbb{R}^3)\cap
W^{1,p}(B,\mathbb{R}^3)$ and $\bar n=n$ in $\Omega$. Mollify $\bar
n$ to get $\bar n^{(j)}\in C^1(\bar\Omega,\mathbb{R}^3)$ with $\bar
n^{(j)}\to n$ in $W^{1,p}(\Omega,\mathbb{R}^3)$. \par Let us define
now an extension of $P$ to the whole $\mathbb{R}^3$, namely the function $\tilde
P:\mathbb{R}^3\to\mathcal{M}^{3\times 3}(\mathbb{R})$ with $\tilde
P(m)\stackrel{\rm{def}}{=}s\left(m\otimes m-\frac{1}{3}Id\right)$. Then

\begin{equation}
\tilde P(\bar n^{(j)})=s\left(\bar n^{(j)}\otimes \bar
n^{(j)}-\frac{1}{3}Id\right)\to \tilde P(n)\textrm{ in }W^{1,p}
\label{eq:napprox}
\end{equation}  since

$$
|\bar n_i^{(j)}\bar n_k^{(j)}-n_in_k|\le |\bar n_{i}^{(j)}(\bar
n_k^{(j)}-n_k)|+|n_k(\bar n_i^{(j)}-n_i)|\le C|\bar n^{(i)}-n|
$$ and, for example,

$$\bar n_i^{(j)}\bar n_{k,\alpha}^{(j)}-n_i n_{k,\alpha}=\underbrace{\bar
n_i^{(j)}(\bar n_{k,\alpha}^{(j)}-n_{k,\alpha})}_{\to
0\,\textrm{ in }L^p}+\underbrace{(\bar n^{(j)}_i-n_i)n_{k,\alpha}}_{\to
0\textrm{ in }L^p}.$$

\par Recalling that  $\bar n^{(j)}$ are $C^1$ functions on $\bar\Omega$ we have:
\begin{equation}
\rm{Tr}\,\tilde P(\bar n^{(j)})=\tilde P(\rm{Tr}\,\bar n^{(j)})
\label{eq:approxtr}
\end{equation}

\par As $\bar n^{(j)}\to n$ in $W^{1,p}(\bar\Omega;\mathbb{R}^3)$ we have that $\rm{Tr}\,\bar n^{(j)}\to\rm{Tr}\,n$ and also that $
\tilde P(\rm{Tr}\,\bar n^{(j)})\to \tilde P(\rm{Tr}\,n)$ in $L^p(\partial\Omega)$. Taking into account this convergence together with (\ref{eq:napprox}), the continuity of the trace operator $\rm{Tr}$ and the fact that $\tilde P(n)=P(n)=Q$ we can pass to the limit (in $L^p(\partial\Omega)$) in both sides of  (\ref{eq:approxtr}) and obtain relation (\ref{eq:boundaryorient}).

\par In order to finish the proof we just need to check that $\rm{Tr}\, n\in L^p(\partial\Omega;\mathbb{S}^2)$. From the definition of the trace we know that $\rm{Tr}\, n\in L^p(\partial\Omega;\mathbb{R}^3)$. Hence recalling (see for instance \cite{evansgariepy}, p. $133$) that we have

\begin{equation}
\lim_{r\to 0}\frac{1}{|B_r(x)\cap \Omega |}\int_{B_r(x)\cap \Omega}|n(y)-\rm{Tr}\,n(y)|\,dy=0, \textrm{ for } \mathcal{H}^2\textrm{  a.e.  }x\in\partial\Omega
\end{equation}  and since $$\int_{B_r(x)\cap\Omega} \Big|1-|\textrm{Tr}\,n|\Big|\,dx=\int_{B_r(x)\cap\Omega} \Big||n|-|\rm{Tr}\,n|\Big|\,dx\le \int_{B_r(x)\cap\Omega} |n-\rm{Tr}\, n|\,dx,$$ we get that $|\rm{Tr}\,n|=1$ $\mathcal{H}^2$ a.e. $x\in\partial\Omega$. $\Box$


\smallskip\begin{remark} One can easily check, by straightforward modifications of the proofs, that  the results of Proposition~\ref{prop:2orient} and Proposition~\ref{boundaryorient}    hold for domains $\Omega\subset\mathbb{R}^2$ with $\boldsymbol{\mathcal{Q}}$ replaced by $\boldsymbol{\mathcal{Q}}_2$.
\label{boundaryorient:2d}
\end{remark}

\par Orientability is preserved by weak convergence:

\begin{proposition} Let $\Omega\subset\mathbb{R}^d$ be a bounded domain with  boundary of class $C^0$.
For $1\le p\le\infty$ let $Q^{(k)}\in
W^{1,p}(\Omega,\boldsymbol{\mathcal{Q}}),k\in\mathbb{N}$ be a
sequence of orientable maps with the corresponding $n^{(k)}\in
W^{1,p}(\Omega,\mathbb{S}^2)$ such that $P(n^{(k)})=Q^{(k)}$. If
$Q^{(k)}$ converges weakly to $Q$ in $W^{1,p}$, where $1\le p<\infty$
(or weak* when $p=\infty$), then $Q$ is orientable and $Q=P(n)$ for
some $n\in W^{1,p}$.\label{limitorient}
\end{proposition}
\smallskip{\bf Proof.} We start by proving an auxiliary result:

\smallskip\begin{lemma}
\label{regularitylemma}
Let $\Omega\subset \mathbb{R}^d$ be an open and bounded set.
If $n\in W^{1,p}(\Omega,\mathbb{S}^2),\,1\le p\le \infty$ then
$Q=P(n)\in W^{1,p}(\Omega,\boldsymbol{\mathcal{Q}})$. Conversely,
assume that $Q\in W^{1,p}(\Omega,\boldsymbol{\mathcal{Q}}),\,1\le
p\le \infty$ and let $n$ be a measurable function on $\Omega$  with
values in $\mathbb{S}^2$ such that $P(n)=Q$. Moreover  assume that $n$
is continuous along almost every line parallel to the coordinate
axes and intersecting $\Omega$. Then $n\in
W^{1,p}(\Omega,\boldsymbol{\mathcal{Q}})$.
\par Moreover:
\begin{equation}
Q_{ij,k}n_j=sn_{i,k}. \label{derivtrick}
\end{equation}
\end{lemma}
 {\bf Proof of the lemma.} For $g,h\in
W^{1,1}(\Omega)\cap L^\infty(\Omega)$ we have $g h\in
W^{1,1}(\Omega)\cap L^\infty(\Omega)$ and
$(gh)_{,i}=gh_{,i}+g_{,i}h$ (see \cite{ADAMS75}). Hence, if $n\in
W^{1,p}(\Omega)$, we have $Q\in W^{1,1}(\Omega)$ and
$Q_{ij,k}=s(n_i n_{j,k}+n_{i,k}n_j)$ from which we obtain
$\nabla Q\in L^p(\Omega)$ and then $Q\in W^{1,p}(\Omega)$. Also
$$Q_{ij,k}n_j=s\left[n_i(n_{j,k}n_j)+n_{i,k}\right]=s[\frac{n_i}{2}(\underbrace{n_jn_j}_{=1})_{,k}
+n_{i,k}]=sn_{i,k}.$$
\par Conversely, suppose that $Q\in W^{1,p}$. Let
$x\in\Omega$ with $n$ continuous along the line
$(x+\mathbb{R}e_k)\cap\Omega$. Let $x+te_k\in\Omega$. As $Q\in
W^{1,1}$ we can suppose that $Q$ is differentiable at $x$ in the
direction $e_k$. Then
\begin{eqnarray}&\hspace{-1in}
\frac{Q_{ij}(x+te_k)-Q_{ij}(x)}{t}=s\left[\cdot\frac{n_i(x+te_k)n_j(x+te_k)-n_i(x)n_j(x)}{t}
\right]\nonumber\\
&=s\cdot n_i(x+te_k)\left[\frac{n_j(x+te_k)-n_j(x)}{t}\right]+
s\cdot\left[\frac{n_i(x+te_k)-n_i(x)}{t}\right]n_j(x).\nonumber
\end{eqnarray}

\par Multiply both sides by
$\frac{1}{2}\left[n_j(x+te_k)+n_j(x)\right]$. Then, since
\begin{eqnarray}&\hspace{-1in}
\left[n_j(x+te_k)-n_j(x)\right]\left[n_j(x+te_k)+n_j(x)\right]\nonumber\\
&=n_j(x+te_k)n_j(x+te_k)-n_j(x)n_j(x)=1-1=0
\end{eqnarray} we have that
\begin{eqnarray}
\frac{Q_{ij}(x+te_k)-Q_{ij}(x)}{t}\cdot
\frac{1}{2}\left[n_j(x+te_k)+n_j(x)\right]&\nonumber\\
&\hspace{-2in}=s\cdot\left[\frac{n_i(x+te_k)-n_i(x)}{t}\right]n_j(x)\frac{1}{2}\left[n_j(x+te_k)+n_j(x)\right].
\end{eqnarray}
Letting $t\to 0$ and using the assumed continuity of $n$ we deduce that
$$s\cdot\lim_{t\to 0}\frac{n_i(x+te_k)-n_i(x)}{t}=Q_{ij,k}(x)n_j(x).$$
\par Hence the partial derivatives of $n$ exist almost everywhere in
$\Omega$ and satisfy
$$sn_{i,k}=Q_{ij,k}n_j$$ and since $\nabla Q\in L^p$ it follows that
$n\in W^{1,p}(\Omega,\mathbb{S}^2)$ as required. $\Box$

\bigskip

\par We continue with the proof of the proposition. As the sequence $Q^{(k)}$
is weakly convergent in $W^{1,p}$ for $1\le p<\infty$ ( weak* for $p=\infty$),
using (\ref{derivtrick}) we have that $n^{(k)}$  is bounded in the
$W^{1,p}$ norm (equi-integrable if $p=1$). Thus there exists
$n\in W^{1,p}(\Omega,\mathbb{S}^2)$ such that on a subsequence
$(n^{(k_l)})_{k_l\in\mathbb{N}}$ we have $n^{(k_l)}\rightharpoonup  n$ in
$W^{1,p}(\Omega)$ and $n^{(k_l)}\to n$ a.e. which implies that
$P(n)=Q$. $\Box$

\smallskip
\begin{remark}  Lemma ~\ref{regularitylemma} shows that in order
to have   $n\in W^{1,p}(\Omega,\mathbb{S}^2)$ such that $P(n)=Q$
for $Q\in W^{1,p}(\Omega,\boldsymbol{\mathcal{Q}})$ it suffices to
have only $n\in W^{1,1}(\Omega,\mathbb{S}^2)$ since then $n$ will have the same regularity as $Q$.
\end{remark}
\medskip
\par The previous result allows us to show that
non-orientability is a stable property with respect to the metric
induced by the distance
$d(Q_1,Q_2)=\|Q_1-Q_2\|_{W^{1,p}(\Omega,\mathbb{R}^9)}$. More
precisely we have:

\begin{lemma}
Let $Q\in W^{1,p}(\Omega,\boldsymbol{\mathcal{Q}}), 1\le p\le\infty$ be
non-orientable. Then there exists $\varepsilon>0$, depending on $Q$,
so that for all $\tilde Q\in W^{1,p}(\Omega,\boldsymbol{\mathcal{Q}})$
with $\|\tilde Q-Q\|_{W^{1,p}(M,\mathbb{R}^9)}<\varepsilon$ the line
field $\tilde Q$ is also non-orientable.
\end{lemma} Taking into account the previous proposition and
reasoning by contradiction the proof is straightforward.

 \bigskip
\par The line-field  theory we have presented  represents
a generalization of the  Oseen-Frank theory, which uses vector fields
for describing uniaxial nematic liquid crystals. The Oseen-Frank theory has been successful in predicting the
equilibrium states as local or global minimizers of an energy
functional:

\begin{equation}
E_{OF}=\int_{\Omega} W(n,\nabla n) \,dx,
\end{equation} where
\begin{eqnarray}
W(n,\nabla n)=K_1(\textrm{div}\,n)^2+K_2(n\cdot\textrm{curl}\,
n)^2+K_3|n\wedge\textrm{curl}\,n|^2\nonumber\\
+(K_2+K_4)(\textrm{tr} (\nabla
n)^2-(\textrm{div}\,n)^2)
\end{eqnarray} and the $K_i$ are elastic constants.
\par   We consider a special case of the Landau~-~de Gennes theory, in which the elastic energy density is
 defined by $$\psi(Q,\nabla Q)=L_1I_1+L_2I_2+L_3I_3+L_4I_4,$$
where  the $L_i$ are constants and the  four elastic invariants
$I_1,\ldots,I_4$ are given by
$$I_1=Q_{ij,j}Q_{ik,k},\,I_2=Q_{ik,j}Q_{ij,k},\,I_3=Q_{ij,k}Q_{ij,k},
\,I_4=Q_{lk}Q_{ij,l}Q_{ij,k},$$
where we have used the summation convention with
$i,j,k\in\{1,2,3\}$.
\par It can be checked that the Oseen-Frank energy is expressible
in terms of the constant $s$ Landau~-~de Gennes $Q$-tensors (see
\cite{newtonmottram}). We have that
$$I_1=s^2\left(|\textrm{div}\,n|^2+|n\wedge\textrm{curl}\,
n|^2\right),\,
I_2=s^2\left(|n\wedge\textrm{curl}\,n|^2+\textrm{tr}\,(\nabla
n)^2\right),$$
$$I_3=2s^2\left(\textrm{tr}\,(\nabla
n)^2+|n\cdot\textrm{curl}\,n|^2+|n\wedge\textrm{curl}\,n|^2\right),$$
$$I_4=2s^3\left(\frac{2}{3}|n\wedge\textrm{curl}\,n|^2-\frac{1}{3}\textrm{tr}(\nabla
n)^2-\frac{1}{3}|n\cdot\textrm{curl}\,n|^2\right).$$

\par We let
$$K_1=L_1s^2+L_2s^2+2L_3s^2-\frac{2}{3}L_4s^3,\;\;K_2=2L_3s^2-\frac{2}{3}L_4s^3,$$
$$K_3=L_1s^2+L_2s^2+2L_3s^2+\frac{4}{3}L_4s^3,\,\;\;K_4=L_2s^2,$$ and observe that the $L_i$ can also be expressed in terms of the $K_i$. Then we have
that
$$W(n,\nabla n)=\psi(Q,\nabla Q),$$ and thus the Oseen-Frank elastic
energy is the same as the Landau~-~de Gennes elastic energy.
 \par For more information concerning the form of the Landau-de
Gennes energy $\psi$ and its relationship to the Oseen-Frank energy
see \cite{dutch},\cite{trebin}.

\section{Orientability issues}
The existence of an oriented,
$\mathbb{S}^2$-valued, version of a
 $\boldsymbol{\mathcal{Q}}$-valued  field is essentially a topological question.
 It amounts to checking if there exists  a lifting of a map that takes
  values into $\boldsymbol{\mathcal{Q}}$ to one
that takes values into its covering space $\mathbb{S}^2$, so that
the lifting has the same regularity as the map.
\par Let us recall (see \cite{hatcher},\cite{lee}) that a map $Q:\Omega\to \boldsymbol{\mathcal{Q}}$ is said to
have a lifting $\varphi^Q:\Omega\to \mathbb{S}^2$ if $P\circ
\varphi^Q=Q$ where $P:\mathbb{S}^2\to \boldsymbol{\mathcal{Q}}$ is a
covering map, defined in our case as
\begin{equation}
P(n)=s(n\otimes n-\frac{1}{3}Id).
\end{equation}
\par In order to show that $P$ is a covering map  we need to  endow
 $\boldsymbol{\mathcal{Q}}$ with an appropriate  topological structure. To this
 end let us first recall that  $\mathbb{R}P^2\stackrel{\rm{def}}{=}\mathbb{S}^2/\sim$ is the
 quotient of $\mathbb{S}^2$ with respect to the equivalence relation
 $n\sim m\Leftrightarrow n=\pm m$ (see \cite{docarmo}, \cite{spivak}).
We define the map $b:\boldsymbol{\mathcal{Q}}\to \mathbb{R}P^2$ by
$$b(s(n\otimes n-\frac{1}{3}Id))=\{n,-n\}\in \mathbb{R}P^2,\,\,\textrm{ for all }n\in\mathbb{S}^2.$$  One can then define in a standard manner
 a topological structure on $\boldsymbol{\mathcal{Q}}$ so that $b$
 is a continuous map. Moreover one can endow
 $\boldsymbol{\mathcal{Q}}$ with a Riemannian structure so that $b$
 is an isometry. Therefore in the remainder of the paper we are able to use for  $Q$-valued
 functions theorems that were proved for
 functions with values into a Riemannian manifold.
\par  The map $P:\mathbb{S}^2\to \boldsymbol{\mathcal{Q}}$ is then
easily seen to be continuous. Also $P$ is surjective and one can
easily check that every point $M\in \boldsymbol{\mathcal{Q}}$ has an
evenly covered neighbourhood (that is there exists an open $U\subset
\boldsymbol{\mathcal{Q}}$, with $M\in U$ and each component of
$P^{-1}(U)$ is mapped homeomorphically onto $U$ by $p$). Thus $P$ is
a covering map (see also \cite{lee} for more details about covering
maps in general).
\par In the case of planar line fields, that is
$\boldsymbol{\mathcal{Q}}_2$-valued fields, one has a similar
lifting problem by identifying, analogously,
$\boldsymbol{\mathcal{Q}}_2$ with $\mathbb{R}P^1$, and denoting
(without loss of generality, regarding $\mathbb{R}P^1$ as embedded
in $\mathbb{R}P^2$ and $\mathbb{S}^1$ in $\mathbb{S}^2$) also with
$P$, the covering map from $\mathbb{S}^1$ to $\mathbb{R}P^1$.
\par There exists a well-developed theory, in algebraic topology, that
shows when it is possible to have a lifting.  Both to avoid the
reader having to enter into the technical details of this theory and
related topological ideas, and for its intrinsic interest, we give
in the next subsection a self-contained treatment of the
orientability of continuous line fields that uses only elementary
point-set topology (but of course with ingredients than have
counterparts in the algebraic topology approach). We use this to
show that for a large class of  two-dimensional domains $G$ one can
check orientability of a continuous line field $Q$ (on $\bar G$)
just by determining the orientability of $Q|_{\partial G}$.
\par However, we face the significant additional difficulty that these topological results are restricted to
continuous functions, while functions in the  Sobolev space
$W^{1,p}(\Omega,\boldsymbol{\mathcal{Q}})$ (with
$\Omega\subset\mathbb{R}^d, \,d=2,3$) are not necessarily continuous
for $p\le d$. Also we are interested only in  liftings that preserve
the Sobolev regularity class. We begin to address these questions in
Section \ref{arbitrary}.

\subsection{Continuous line fields on arbitrary domains}
\label{subsection:alternative}
\par We first show how to lift continuous line fields along continuous paths:

\begin{lemma} If $-\infty<t_1<t_2<\infty$ and
$Q:[t_1,t_2]\to\boldsymbol{\mathcal{Q}}$ is continuous then there
exist exactly two continuous maps $($liftings$)$
$n^+,n^-:[t_1,t_2]\to\mathbb{S}^2,\,$ so that
\begin{equation}Q(t)=s(n^\pm(t)\otimes n^\pm(t)-\frac{1}{3}Id),
\end{equation} and $n^+=-n^-$. $($Equivalently, given either of the two possible initial orientations $\bar n\in \mathbb{S}^2$, so that $Q(t_1)=s(\bar n\otimes \bar n-\frac{1}{3}Id)$, there exists a unique continuous lifting $n:[t_1,t_2]\to\mathbb{S}^2$ with $n(t_1)=\bar n$.$)$
\par Suppose in addition that $\bar Q=s(\bar m\otimes \bar
m-\frac{1}{3}Id),\bar m\in\mathbb{S}^2$, and that $|Q(t)-\bar Q|\le
\varepsilon |s|$ for all $t\in [t_1,t_2]$ where
$0<\varepsilon<\sqrt{2}$. Then one of the liftings, $n^+$ say,
satisfies
\begin{equation}\label{alternative0}
|n^+(t)-\bar m|\le \varepsilon,\,\mbox{ for all } t\in [t_1,t_2]
\end{equation} and the other $n^-=-n^+$  satisfies

\begin{equation}\label{alternative0a}
|n^-(t)+\bar m|\le \varepsilon,\,\mbox{ for all } t\in [t_1,t_2]
\end{equation}
\label{lemma:alternative1}
\end{lemma}
\smallskip{\bf Proof.} Let $0<\varepsilon<\sqrt{2}$. Given
$n,\bar m\in\mathbb{S}^2$ with $|n\otimes n-\bar m\otimes \bar
m|\le\varepsilon$ we have that $|n\otimes n-\bar m\otimes \bar
m|^2=2(1-(n\cdot \bar m)^2)\le\varepsilon^2$ and so
\begin{equation}n\cdot\bar m\ge \sqrt{1-\frac{\varepsilon^2}{2}}>0\textrm{ or }n\cdot \bar m\le
-\sqrt{1-\frac{\varepsilon^2}{2}}<0. \label{alternative4}
\end{equation}

\noindent Thus $n\otimes n=n^+\otimes n^+=n^-\otimes n^-$, where
$n^+\cdot\bar m>0$ and $n^-=-n^+$ satisfies $n^-\cdot \bar m<0$.
\par Now let $Q(\tau)=s(n(\tau)\otimes n(\tau)-\frac{1}{3}Id)$ be continuous on $[t_1,t_2]$. Then
there exists $\delta>0$ such that $|n(\tau)\otimes
n(\tau)-n(\sigma)\otimes n(\sigma)|<\sqrt{2}$ for all
$\sigma,\tau\in [t_1,t_2]$ with $|\sigma-\tau|\le\delta$, and we may
suppose that $t_2-t_1=M\delta$ for some integer $M$. First take
$\bar m\stackrel{\rm def}{=}n(t_1)$ and for each $\tau\in
[t_1,t_1+\delta]$ choose $n^+(\tau)$ as above so that
$n^+(\tau)\otimes n^+(\tau)=n(\tau)\otimes n(\tau)$ and
$n^+(\tau)\cdot\bar m>0$. We claim that
$n^+:[t_1,t_1+\delta]\to\mathbb{S}^2$ is continuous. Indeed let
$\sigma_j \to \sigma$ in $[t_1,t_1+\delta]$ and suppose for
contradiction that $n^+(\sigma_j)\not\to n^+(\sigma)$. Then since
$n^+(\sigma_j)\otimes n^+(\sigma_j)\to n^+(\sigma)\otimes
n^+(\sigma)$ there is a subsequence $\sigma_{j_k}$ such that
$n^+(\sigma_{j_k})\to -n^+(\sigma)$. But then $-n^+(\sigma)\cdot
m\ge 0$ a contradiction which proves the claim. Repeating this
procedure with $\bar n=n^+(t_1+\delta)$ we obtain a continuous
lifting $n^+:[t_1,t_1+2\delta]\to\mathbb{S}^2$, and thus inductively
a continuous lifting $n^+:[t_1,t_2]\to\mathbb{S}^2$. Setting
$n^-=-n^+$ gives a second continuous lifting.
\par If $n^*:[t_1,t_2]\to\mathbb{S}^2$ is a continuous lifting then we may suppose that $n^*(t_1)=n^+(t_1)$ say.
 We claim that then $n^*(\tau)=n^+(\tau)$ for all $\tau\in [t_1,t_2]$. If
 not, by continuity there would be a first $T>t_1$ with
 $n^*(T)=n^-(T)$. But then $n^*(T)=\lim_{
\tau\to T_-}n^*(\tau)=\lim_{\tau\to
 T_-}n^+(\tau)=n^+(T)$, a contradiction. Thus there are exactly two
 continuous liftings.
 \par Finally suppose that $|Q(t)-\bar Q|\le \varepsilon |s|$ for
 all $t\in [t_1,t_2]$ with $0<\varepsilon<\sqrt{2}$. Then by
 (\ref{alternative4})  and the continuity of the lifting one of the liftings, $n^+$ say, satisfies
 $n^+(t)\cdot \bar m\ge \sqrt{1-\frac{\varepsilon^2}{2}}$ and so
 $|n^+(t)-\bar m|^2=2(1-n^+(t)\cdot\bar m)\le
 2(1-\sqrt{1-\frac{\varepsilon^2}{2}})\le\varepsilon^2$ and the
 result follows. $\Box$\
 \begin{proposition}
\label{continuouslift}
Let $\Omega\subset\mathbb{R}^d$ be a simply-connected domain, and let $Q:\Omega\to\mathcal{Q}$ be continuous. Then there exists a continuous lifting $n_Q:\Omega\to\mathbb{S}^2$.
\end{proposition}
{\bf Proof.} This is standard, and follows, for example,
 from \cite{hatcher} p.~61, Prop.~1.33. To give a direct argument,
  fix $x_0\in \Omega$ and let $n^0\in \mathbb{S}^2$ be one of the two possible orientations
  of $Q(x_0)$. Given any $x\in \Omega$, let $\gamma:[0,1]\to\Omega$ be a continuous path
  with $\gamma(0)=x_0, \gamma(1)=x$. By Lemma~\ref{lemma:alternative1} there is a unique
   continuous lifting $n:[0,1]\to\mathbb{S}^2$ of $Q(\gamma(\cdot))$ with $n(0)=n^0$.
    Define $n^Q(x)=n(1)$. Then the method of proof of Theorem~\ref{theorem:alternative} below,
     which for simply-connected domains does not depend on the dimension $d$, shows that $n^Q$
     is well defined and continuous.  (Note that local path connectedness holds because $\Omega$ is open,
      and that none of the complications concerning the $\omega_i$ are needed, so
      that $h^*(\lambda,\cdot)=h(\lambda,\cdot)$.)  $\Box$.

\smallskip\par In the rest of this subsection we restrict ourselves to a class of topologically non-trivial domains in $2D$ and
 show that one can check orientability at the boundary.
\smallskip\par  Let $\Omega\subset\mathbb{R}^2$ be a
bounded  domain with $\partial\Omega$ a Jordan curve. For $1\le i\le
N$ let $\omega_i\subset\mathbb{R}^2$ be a bounded   domain with
$\partial\omega_i$ a Jordan curve and
$\overline{\omega_i}\subset\Omega,\,\overline{\omega_i}\cap\overline{\omega_j}=\emptyset$
if $i\not=j$. Note that $\Omega$ and $\omega_i$ are simply connected
(for example by the proof of Lemma~\ref{lemma:alternative3}~(iii)).
Let
\begin{equation}
\label{holes}G=\Omega\backslash\cup_{i=1}^N \overline{\omega_i}.
\end{equation}
  We call such a domain $G$ a {\it domain with holes}.

\smallskip\begin{theorem} Let $Q:\overline
G\to\boldsymbol{\mathcal{Q}}_2$ be continuous with
$Q|_{\partial\omega_i}$ orientable as a continuous function for
$1\le i\le N$. Then $Q$ is orientable as a continuous function.
\label{theorem:alternative}
\end{theorem}

\par The proof of Theorem~\ref{theorem:alternative} requires some preparation. Let $f:[0,1]\to\overline{\Omega}$ be a continuous path with
$f(0)$ and $f(1)$ belonging to $\overline{G}$. We define a new
{continuous} path $f^*:[0,1]\to\overline{G}$ with $f^*(0)=f(0),
f^*(1)=f(1)$ by replacing the parts of $f$ where it lies in some
$\omega_i$ by corresponding paths on $\partial \omega_i$. If
$f(t)\in\overline{G}$ we set $f^*(t)=f(t)$. Otherwise consider an
interval $[t_1,t_2]\subset [0,1]$ such that $f(t)\in\omega_i$ for
$t_1<t<t_2,
\,f(t_1)\in\partial\omega_i,\,f(t_2)\in\partial\omega_i$, for some
$i\in\{1,\dots,N\}$.
\par Let $\gamma_i:\mathbb{S}^1\to\mathbb{R}^2$ parametrize
$\partial\omega_i$, so that $\gamma_i=\gamma_i(\theta)$ can be
identified with a $2\pi$-periodic function
$\gamma_i:\mathbb{R}\to\mathbb{R}^2$,
$\gamma_i(\theta+2\pi)=\gamma_i(\theta)$, where $\theta$ denotes the
polar angle of $\mathbb{S}^1$. Then
$\gamma_i(\theta_1)=f(t_1),\gamma_i(\theta_2)=f(t_2)$, where
$\theta_i\in [0,2\pi)$. For $t\in (t_1,t_2)$ we define
\begin{equation}
f^*(t)=\gamma_i(\theta(t)) \label{alternative5}
\end{equation} where
\begin{equation}
\theta(t)=\theta_1+(\tilde\theta_2-\theta_1)\left(\frac{t-t_1}{t_2-t_1}\right)
\label{alternative6}
\end{equation} and
\begin{equation}\label{alternative6a}
\tilde\theta_2=\left\{\begin{array}{ll} \theta_2
&\textrm{ if }|\theta_2-\theta_1|\le\pi\\
\theta_2-\textrm{sgn}(\theta_2-\theta_1)2\pi &\textrm{ if
}|\theta_2-\theta_1|>\pi\end{array}\right. . \end{equation}

\par Thus $f^*$ is continuous from $[t_1,t_2]\to\partial\omega_i$
and traces the image under $\gamma_i$ of the minor (shorter) arc
joining $\theta_1$ and $\theta_2$ (with an unimportant specific
choice of the arc if $\theta_1,\theta_2$ represent opposite points
of $\mathbb{S}^1$).

\smallskip\begin{lemma}
$f^*:[0,1]\to\overline{G}$ is continuous. \label{lemma:alternative2}
\end{lemma}
{\bf Proof.} Let $\sigma_j\to\sigma$ in $[0,1]$. We show that
$f^*(\sigma_j)\to f^*(\sigma)$. If $f(\sigma)\in
G\cup\partial\Omega$ then this is obvious, since then
$f(\sigma_j)\in\overline{G}$ for sufficiently large $j$, and so
\begin{equation}
\lim_{j\to\infty}f^*(\sigma_j)=\lim_{j\to\infty}f(\sigma_j)=f(\sigma)=f^*(\sigma).
\label{alternative8}
\end{equation}
\par If $f(\sigma)\in\omega_i$ then, since $f(0)=f(1)\in\bar G$, there exists an interval $[t_1,t_2]\subset[0,1]$ as above with $t_1<\sigma<t_2$, so that $\lim_{j\to\infty}f^*(\sigma_j)=f^*(\sigma)$ follows from \eqref{alternative5}-\eqref{alternative6a}. Thus, arguing by contradiction, we may assume that $f(\sigma)\in\partial\omega_i$ for some $i$.
If $f(\sigma_j)\in\partial\omega_i$ for all $j$ then again we have
(\ref{alternative8}). Thus it suffices to consider the case when
$f(\sigma_j)\in\omega_i$ for all $j$ and either $(i)
\sigma_j>\sigma$ for all $j$, or $(ii) \sigma_j<\sigma$ for all $j$.
\par We assume $(i)$ with $(ii)$ being treated similarly. If
$f(t)\in\omega_i$ for all $t>\sigma$ with $t-\sigma$ sufficiently
small, then we can take $t_1=\sigma$ and again   $\lim_{j\to\infty}f^*(\sigma_j)=f^*(\sigma)$ follows from \eqref{alternative5}-\eqref{alternative6a}.
\par Otherwise there exist $t_{2j}>\sigma_j>t_{1j}>\sigma$ with
$f(t)\in\omega_i$ for all $t\in (t_{1j},t_{2j}),
f(t_{1j})\in\partial\omega_i,\,f(t_{2j})\in\partial\omega_i$ and
$t_{2j}\to \sigma$. We have that
$f(t_{1j})=\gamma_i(\theta_{1j}), f(t_{2j})=\gamma_i(\theta_{2j})$ with
$\theta_{1j},\theta_{2j}\in [0,2\pi)$ and
$\theta_{1j}\stackrel{\mathbb{S}^1}{\to}\theta,
\theta_{2j}\stackrel{\mathbb{S}^1}{\to}\theta$, where
$\gamma_i(\theta)=f(\sigma)$, since otherwise, for example, we would
have a subsequence $\theta_{1j_k}$ with
$\theta_{1j_k}\stackrel{\mathbb{S}^1}{\to}\bar\theta\not=\theta \,\mbox{mod}\,(2\pi)$ and
$\gamma_i(\theta)=\gamma_i(\bar\theta)$, contradicting that
$\partial\omega_i$ is a Jordan curve. Thus from
(\ref{alternative5}), (\ref{alternative6}) we have
$$f^*(\sigma_j)=\gamma_i\left(\theta_{1j}+(\tilde\theta_{2j}-\theta_{1j})
\left(\frac{\sigma_j-t_{1j}}{t_{2j}-t_{1j}} \right)\right)$$ and so
$\lim_{j\to\infty}
f^*(\sigma_j)=\gamma_i(\theta)=f(\sigma)=f^*(\sigma)$ as required.
$\Box$

\smallskip\begin{lemma}\hspace{1in} \newline \vspace{-.15in}
\par{$($i$)$} $\overline{\Omega}$ is path-connected and simply connected.
\par{$($ii$)$} $\overline{G}$ is path-connected.
\par$($iii$)$ $\bar G$ is locally path-connected: given any $x\in\overline{G}$ and $\bar\varepsilon>0$,
there exists $\delta>0$ such that for any $z\in\overline{G}$ with
$|x-z|<\delta$ there is a  continuous path
$\tilde\gamma:[0,1]\to\overline{G}$ with
$\tilde\gamma(0)=x,\tilde\gamma(1)=z$ and
$|\tilde\gamma(t)-x|<\bar\varepsilon$ for all $t\in [0,1]$.
\label{lemma:alternative3}
\end{lemma}
\smallskip{\bf Proof.} By the Schoenflies theorem \cite{bing}
there is a homeomorphism $u:\mathbb{R}^2\to\mathbb{R}^2$ with
$u(D)=\Omega, u(\partial D)=\partial\Omega$ and
$u(\overline{D}^c)=\overline{\Omega}^c$ where $D=B(0,1)$ is the unit
disk. Since $\overline{D}$ is path connected and simply connected
this proves $(i)$. Given $x,z\in\overline{G}$ there thus exists a
continuous path $f:[0,1]\to\overline{\Omega}$ with $f(0)=x,f(1)=z.$
By Lemma~\ref{lemma:alternative2}, $f^*:[0,1]\to\overline{G}$ is
continuous with $f^*(0)=x,f^*(1)=z$, which proves $(ii)$. To prove
$(iii)$, note that this is obvious if $x\in G$ so we may suppose that
$x\in\partial\Omega$. The argument in the case
$x\in\partial\omega_i$ for some $i$ is  similar. Note that
$u^{-1}(x)\in\partial D$. Let $\sigma>0$ be sufficiently small so
that $B(u^{-1}(x),\sigma)\cap \cup_{i=1}^N u^{-1}(\partial\omega_i)$
is empty and $|y-u^{-1}(x)|<\sigma$ implies
$|u(y)-x|<\bar\varepsilon$. Then let $\delta>0$ be such that
$|z-x|<\delta$ implies $|u^{-1}(z)-u^{-1}(x)|<\sigma$. Then
$\tilde\gamma(t)=u(tu^{-1}(z)+(1-t)u^{-1}(x))$ defines a suitable
path.  $\Box$

\smallskip
\par {Before starting the proof of
Theorem~\ref{theorem:alternative} we need one more technical lemma:
\begin{lemma}
\label{lemma:alternative1a} Let $G$ be a domain with holes as above
and let $Q\in C(\bar G,\boldsymbol{\mathcal{Q}})$. There exists $\nu>0$
such that if $\bar x\in \bar G$, $-\infty<t_1<t_2<\infty$,
$f^{(j)},f\in C([t_1,t_2];\bar G)$, $f^{(j)}([t_1,t_2])\subset
B(\bar x,\nu)$, $f^{(j)}(t_2)\to f(t_2)$, and if
$n^{(j)},n:[t_1,t_2]\to\mathbb{S}^2$ are continuous liftings of
$Q(f^{(j)}(\cdot)), Q(f(\cdot))$ respectively with
$|n^{(j)}(t_1)-n(t_1)|<1$, then $n^{(j)}(t_2)\to n(t_2)$.
\end{lemma}
\smallskip{\bf Proof.}
Choose $\nu$ sufficiently small such that $|Q(x)-Q(y)|\leq |s|$ if
$x,y\in\bar G$ {with $|x-y|\leq 2\nu$}. In
Lemma~\ref{lemma:alternative1} set $\bar m=n(t_1)$. Then by
\eqref{alternative0}, \eqref{alternative0a} we have that
\begin{equation}\label{alternative0b}
|n(t)+n(t_1)|>1  \mbox{ for all } t\in[t_1,t_2].
\end{equation}
Also
\begin{equation}
|n^{(j)}(t_1)+n(t_1)|\geq 2|n(t_1)|-|n^{(j)}(t_1)-n(t_1)|>1
\end{equation}
 and so by Lemma~\ref{lemma:alternative1}
\begin{equation}\label{alternative0c}
|n^{(j)}(t)-n(t_1)|\leq 1 \mbox{ for all }t\in[t_1,t_2].
\end{equation}
Suppose that $n^{(j)}(t_2)\not\to n(t_2)$. Then since
$Q(f^{(j)}(t_2))\to Q(f(t_2))$ there exists a subsequence $j_k$ such
that
 $n^{(j_k)}(t_2)\to - n(t_2)$.
 But then, from \eqref{alternative0c},
$|n(t_2)+n(t_1)|\leq 1$, contradicting \eqref{alternative0b}. $\Box$

\smallskip
\noindent{\bf Proof of Theorem~\ref{theorem:alternative}.} Let
$x_0\in G$ and choose one of the two possible orientations $(m^0,0)$
for $Q(x_0)$, where $m^0\in\mathbb{S}^1$. Let $x\in\overline{G}$ be
arbitrary. By Lemma~\ref{lemma:alternative3}$(ii)$ there
exists a continuous path $\gamma:[0,1]\to\overline G$ with
$\gamma(0)=x_0, \,\gamma(1)=x$. By Lemma~\ref{lemma:alternative1}
there exists a unique continuous lifting $n:[0,1]\to\mathbb{S}^1$
such that $Q(\gamma(t))=s\left((n(t),0)\otimes
(n(t),0)-\frac{1}{3}Id\right),\,t\in [0,1]$ and $n(0)=m^0$. We
define $N(x)\stackrel{\rm def}{=}n(1)$. To show that $N(x)$ is well defined, suppose
that $\gamma':[0,1]\to\overline{G}$ is another continuous path with
$\gamma'(0)=x_0,\gamma'(1)=x$, let $n':[0,1]\to\mathbb{S}^1$ be the
corresponding continuous lifting, and suppose for contradiction that
$n'(1)\not= N(x)$, so that $n'(1)=-N(x)$. Define the continuous loop
$\Gamma:[0,1]\to\overline{G}$ by

$$\Gamma(t)=\left\{\begin{array}{ll}\gamma(2t) & \textrm{ if }0\le
t\le\frac{1}{2}\\
\gamma'(2(1-t)) &\textrm{ if }\frac{1}{2}\le t\le
1\end{array}\right.$$ so that $\Gamma(0)=\Gamma(1)=x_0.$ Then
$\tilde N(t)$ defined by

$$\tilde N(t)=\left\{\begin{array}{ll} n(2t) & \textrm{ if }0\le
t\le\frac{1}{2}\\
-n'(2(1-t)) &\textrm{ if }\frac{1}{2}\le t\le 1\end{array}\right.$$
is a continuous lifting such that

\begin{equation}
Q(\Gamma(t))=s\left((\tilde N(t),0)\otimes (\tilde
N(t),0)-\frac{1}{3}Id\right) \label{alternative9} \end{equation}
and $\tilde N(0)=m^0,\tilde N(1)=-m^0$. The bulk of the proof will
be to show that such a lifting cannot exist, so that $N(x)$ is well
defined. Assuming this, we claim that
$N:\overline{G}\to\mathbb{S}^1$ is a continuous lifting of $Q$, that
is
$$Q(x)=s\left((N(x),0)\otimes (N(x),0)-\frac{1}{3}Id\right)$$ for
all $x\in\overline{G}$ and $N$ is continuous in $x$. We only need to
prove the continuity.

Let $z_j\in\overline{G}$ with $z_j\to x$. By Lemma~\ref{lemma:alternative3}~(iii) there exist continuous
 paths  $\tilde\gamma_j:[0,1]\to\bar G$ with {$\tilde\gamma_j(0)=x, \tilde\gamma_j(1)=z_j$} and $\tilde\gamma_j(\cdot)\to x$ in $C([0,1];\bar G)$. Then we may consider the path $$\hat\gamma_j(t)=\left\{\begin{array}{ll} \gamma(2t) &\textrm{ if
}0\le t\le\frac{1}{2}\\
\tilde\gamma_j(2t-1) &\textrm{ if }\frac{1}{2}\le t\le
1\end{array}\right.$$  so that $N(x)=\tilde n^j(0), N(z_j)=\tilde
n^j(1)$, where $\tilde n^j:[0,1]\to\mathbb{S}^1$ is such that
$(\tilde n^j,0)$ is the unique continuous lifting of
$Q(\tilde\gamma_j(\cdot))$ with $\tilde n^j(0)=N(x)$.  By
Lemma~\ref{lemma:alternative1a} we deduce that $N(z_j)\to N(x)$ as
required.

\par It remains to prove that there is no continuous loop
$\Gamma:[0,1]\to\overline{G}$ with $\Gamma(0)=\Gamma(1)=x_0$ and a
corresponding continuous lifting $\tilde N:[0,1]\to\mathbb{S}^1$ so
that (\ref{alternative9}) holds and $\tilde N(0)=m^0,\tilde
N(1)=-m^0$.
\par Since $\overline{\Omega}$ is simply connected, there exists a
continuous homotopy $h:[0,1]^2\to\overline{\Omega}, h=h(\lambda,t)$
with $h(0,t)=\Gamma(t),h(1,t)=x_0$ for all $t\in [0,1]$ and
$h(\lambda,0)=h(\lambda,1)=x_0$ for all $\lambda\in [0,1]$. For each
$\lambda$ we consider the path
$h^*(\lambda,\cdot):[0,1]\to\overline{G}$, which is continuous by
Lemma~\ref{lemma:alternative2}. By Lemma~\ref{lemma:alternative1},
for each $\lambda\in [0,1]$ there is a unique continuous lifting
$n^{\lambda}:[0,1]\to\mathbb{S}^1$ such that

\begin{equation}
Q(h^*(\lambda,t))=s\left((n^{\lambda}(t),0)\otimes
(n^{\lambda}(t),0)-\frac{1}{3}Id\right)
\end{equation} and $n^\lambda(0)=m^0$. We know that $n^0(1)=\tilde
N(1)=-m^0$. We will prove that $n^\lambda(1)$ is a continuous
function of $\lambda\in [0,1]$ so that $n^\lambda(1)=-m^0$ for all
$\lambda\in [0,1]$. In particular $n^1(1)=-m^0$ contradicting
$h(1,t)=x_0$ for all $t$ and $n^1(0)=m^0$ (since $n^1(t)\in \{m^0,
-m^0\}$ is continuous). This contradiction proves the theorem.
\par To prove the continuity of $n^\lambda(1)$ in $\lambda$ we make
use of the assumption that $Q|_{\partial\omega_i}$ is orientable for
each $i$. Let $\lambda_k\to\lambda$ in $[0,1]$. Define
$$T=\sup \{t\in [0,1]: h(\lambda,t)\in G\cup\partial\Omega,
n^{\lambda_k}(t)\to n^\lambda(t)\}.$$
 Since $h(\lambda,t)\in\overline{G}$ for all $\lambda\in [0,1]$
and all sufficiently small $t$ it follows from
Lemma~\ref{lemma:alternative1a} that $T>0$. Suppose for contradiction
that $T<1$. If $h(\lambda,T)\in G\cup\partial\Omega$ then by the
continuity of $h$ there exists a $\sigma>0$ such that
$h(\lambda_k,t)\in B(h(\lambda,T),\delta)$ for all sufficiently
large $k$ and for all $t$ with $|t-T|<\sigma$, where $\delta>0$ is
small enough so that
$B(h(\lambda,T),\delta)\cap\overline{\Omega}\subset\overline{G}$ and $\delta
<\nu$, where $\nu$ is as given in Lemma~\ref{lemma:alternative1a}.
  Hence if $|t-T|<\sigma$,
$h(\lambda,t)=\lim_{k\to\infty} h(\lambda_k,t)\in\overline{G}$, and
so $h^*(\lambda,t)=h(\lambda,t)$. Also by the definition of $T$,
there exists $\tau\in (T-\sigma,T)$ with $n^{\lambda_k}(\tau)\to
n^\lambda(\tau)$, and so by Lemma~\ref{lemma:alternative1a}   $n^{\lambda_k}(t)\to n^\lambda(t)$ for $t\in(T,T+\sigma)$,
contradicting the definition of $T$.
\par Thus we may suppose that $h(\lambda,T)\in \partial\omega_i$ for
some $i$. Let $[T,\bar T]$ be the maximal closed interval,
containing $T$, such that $h(\lambda,t)\in\overline{\omega_i}$ for
all $t\in [T,\bar T]$.

\smallskip\begin{lemma} $n^{\lambda_k}(T)\to
n^\lambda(T),\,n^{\lambda_k}(\bar T)\to n^\lambda(\bar T)$
\end{lemma}
\smallskip{\bf Proof.}
{\it Step 1.} We show that given $\delta>0$ there exist
$\sigma>0$ and $k_0$ such that
\begin{equation}
|h^*(\lambda_k,t)-h(\lambda,T)|<\delta \label{alternative12}
\end{equation} whenever $|t-T|<\sigma$ and $k\ge k_0$.
\par If this were not true there would be a sequence $t_j\to T$ and
a subsequence $k_j\to \infty$ such that
$$|h^*(\lambda_{k_j},t_j)-h(\lambda,T)|\ge \delta\,\mbox{ for all } j.$$
 Since $\lim_{j\to\infty} h(\lambda_{k_j},t_j)=h(\lambda,T)$ we
may suppose that $h(\lambda_{k_j},t_j)\in\omega_i$ for all $j$, and thus we have $t_{1j}<t_j<t_{2j}$ with
$h(\lambda_{k_j},t)\in\omega_i$ for $t\in (t_{1j},t_{2j})$,
$h(\lambda_{k_j},t_{1j}),
h(\lambda_{k_j},t_{2j})\in\partial\omega_i$. By the definition of
$T$ there exists a sequence $T_l\to T-$ such that $h(\lambda, T_l)\in
G$ (and $n^{\lambda_k}(T_l)\to n^{\lambda}(T_l)$). In particular
$h(\lambda_{k_j},T_l)\in G$ for $j$ sufficiently large, and so
$t_{1j}\ge T_l$ for large enough $j$. Hence $t_{1j}\to T$ and so
$h^*(\lambda_{k_j}, t_{1j})=h(\lambda_{k_j},t_{1j})\to h(\lambda, T)$. Thus $h^*(\lambda_{k_j},t_{1j})=\gamma_i(\theta_{1j}), h^*(\lambda_{k_j},t_{2j})=\gamma_i(\theta_{2j})$, where $\theta_{1j}\to\theta$ and $h(\lambda,T)=\gamma_i(\theta)$. From \eqref{alternative5}-\eqref{alternative6a}
$$h^*(\lambda_{k_j},t_j)=\gamma_i\left(\theta_{1j}+(\tilde\theta_{2j}-\theta_{1j})
)\left(\frac{t_j-t_{1j}}{t_{2j}-t_{1j}}\right)\right).$$
 Considering separately the cases when
$t_{2j}-t_{1j}\to 0$ and when $t_{2j}-t_{1j}\ge \mu>0$ we get
$\lim_{j\to\infty} h^*(\lambda_{k_j},t_j)=h(\lambda,T)$.
\par This contradiction proves the claim, which by a similar
argument also holds if $T$ is replaced by $\bar T$.\vspace{.1in}

\noindent{\it Step $2$.} We prove that $n^{\lambda_k}(T)\to
n^{\lambda}(T)$. In Step $1$, we choose $\delta\in(0,\nu)$, where $\nu$ is given in Lemma~\ref{lemma:alternative1a}, and note that $h^*(\lambda_k,T)\to h^*(\lambda,T)=h(\lambda,T)$.
  Since  $n^{\lambda_k}(T_l)\to n^{\lambda}(T_l)$ the result follows from Lemma~\ref{lemma:alternative1a} applied on the interval $[T_l,T]$.  \vspace{.1in}

 \noindent {\it Step $3$.} We prove that $n^{\lambda_k}(\bar T)\to
n^{\lambda}(\bar T)$.
\par If $T=\bar T$ there is nothing to prove and so we assume $\bar
T>T$. First we note that
\begin{equation}
\lim_{k\to\infty} \sup_{t\in [T,\bar T]}
\textrm{dist}(h^*(\lambda_k,t),\partial\omega_i)=0,
\label{alternative13}
\end{equation}
since otherwise there would exist a subsequence $k_j\to\infty$
and $t_j\to t$ in $[T,\bar T]$ with
$h^*(\lambda_{k_j},t_j)=h(\lambda_{k_j},t_j)\to
h(\lambda,t)\not\in\overline{\omega_i}$, a contradiction.

\par Since $Q|_{\partial\omega_i}$ is orientable, there is a unique
continuous lifting $\hat N:\partial\omega_i\to\mathbb{S}^1$ such
that $\hat N(h(\lambda,T))=n^\lambda(T)$. Given $\varepsilon>0$
sufficiently small, we choose $\delta\in(0,\varepsilon), k_0$ such
that (from (\ref{alternative13})) if $k\geq   k_0$
\begin{equation}
\textrm{dist}(h^*(\lambda_k,t),\partial\omega_i)<\delta,\,\mbox{ for all }
t\in[T,\bar T], \label{alternative14}
\end{equation}
\begin{equation}
x,y\in\bar G\textrm{ with }|x-y|<4\delta\textrm{ implies
 }|Q(x)-Q(y)|<\varepsilon|s|,
 \label{alternative15}
\end{equation}
\begin{equation}
|\hat N(z)-\hat N(\bar z)|<\varepsilon, \textrm{ if }z,\bar
z\in\partial\omega_i\textrm{ with }|z-\bar z|<3\delta,
\label{alternative16}
\end{equation} and
\begin{equation}
|n^{\lambda_k}(T)-n^{\lambda}(T)|<\varepsilon,\,|h^*(\lambda_k,T)-h(\lambda,T)|\le\delta.
\label{alternative17}
\end{equation}
 For $k\geq k_0$ define
\begin{eqnarray*}S_k=\{t\in [T,\bar T]:\textrm{ there exists
}z=z(t)\in\partial\omega_i &\\& \hspace{-2in}\textrm{ with }
|z-h^*(\lambda_k,t)|\le\delta, |n^{\lambda_k}(t)-\hat N(z)|\le
2\varepsilon\}.
\end{eqnarray*}
 It is easily seen that $S_k$ is closed. Also $T\in S_k$ because we
can take $z=h(\lambda,T)$ and use (\ref{alternative17}). We show
that if $t\in S_k$ with $t<\bar T$ then $t+s\in S_k$ for $s>0$
sufficiently small, so that $S_k=[T,\bar T]$.
\par Given $s$, by (\ref{alternative14}) there exists
$z(t+s)\in\partial\omega_i$ with
$$|z(t+s)-h^*(\lambda_k,t+s)|<\delta.$$
 If $s$ is chosen small enough so that
$|h^*(\lambda_k,t+s)-h^*(\lambda_k,t)|<\delta$ we have that
\begin{equation}
|h^*(\lambda_k,t+s)-z(t)|\le 2\delta,
\end{equation}
\begin{equation}
|z(t+s)-z(t)|< 3\delta. \label{alternative19}
\end{equation}
 Thus by (\ref{alternative15})
$$|Q(h^*(\lambda_k,t+s))-Q(z(t))|<\varepsilon|s|.$$
 Also $|n^{\lambda_k}(t)-\hat N(z(t))|\le 2\varepsilon$ and so by
Lemma~\ref{lemma:alternative1} we have $|n^{\lambda_k}(t+s)-\hat
N(z(t))|\le \varepsilon$ and hence, by
(\ref{alternative16}),(\ref{alternative19})
$|n^{\lambda_k}(t+s)-\hat N(z(t+s))|\le 2\varepsilon$. Hence $t+s\in
S_k$ as required.
\par Since $\bar T\in S$ and since $h^*(\lambda_k,\bar T)\to
h(\lambda,\bar T)$ letting $\varepsilon\to 0$ we deduce that
$n^{\lambda_k}(\bar T)\to \hat N(h(\lambda,\bar T))$. But
$h^*(\lambda,t)\in\partial\omega_i$ and is continuous in $t$ for
$t\in [T,\bar T]$. Hence $n^{\lambda}(t)=\hat N(h^*(\lambda,t))$ for
all $t\in [T,\bar T]$ and in particular $\hat N(h(\lambda,\bar
T))=n^{\lambda}(\bar T)$. Thus $n^{\lambda_k}(\bar T)\to
n^{\lambda}(\bar T)$ as required. $\Box$

\smallskip\par To complete the proof of the theorem we note that by
the definition of $\bar T$ there exists a sequence $\bar T_r\to \bar
T+$  with $h(\lambda,\bar T_r)\in G$. By (\ref{alternative12}) for
$\bar T$, given $\delta\in(0,\nu)$, $\nu$ as in
Lemma~\ref{lemma:alternative1a},  we have
$|h^*(\lambda_k,t)-h(\lambda,\bar T)|<\delta$ for $k\ge k_0$ and
$|t-\bar T|<\sigma$. Let $r$ be large enough so that $|\bar T_r-\bar
T|<\sigma$.  Then for $k$ large enough $h^*(\lambda_k,\bar
T_r)=h(\lambda_k,\bar T_r)\in G$, and so $h^*(\lambda_k,\bar T_r)\to
h^*(\lambda,\bar T_r)$.  Since also $n^{\lambda_k}(\bar T)\to
n^{\lambda}(\bar T)$ we deduce from Lemma~\ref{lemma:alternative1a}
that  $n^{\lambda_k}(\bar T_r)\to n^{\lambda}(\bar T_r)$ for $r$
sufficiently large, contradicting the definition of $T$. Thus $T=1$,
and using the same argument as just after the definition of $T$ we
deduce that $n^{\lambda_k}(1)\to n^{\lambda}(1)$. $\Box$

\begin{remark}
Theorem~\ref{theorem:alternative} could alternatively have been proved using algebraic topology notions. The orientability
of a continuous line field on $\bar G$ needs to be checked only on a
set of generators of the fundamental group $\pi_1(\bar G)$ of $\bar G$. It seems to
be well known (though tricky to find  in the literature) that for a
domain with holes as defined before, the boundary loops
$\partial\omega_i, i=1,\dots,n$ constitute a family of generators of
$\pi_1(\bar G)$. These observations suffice to give the result of
Theorem~\ref{theorem:alternative}.
\end{remark}

\subsection{Arbitrary line fields on  simply-connected domains}
\label{arbitrary}

 In this case we have that there is a lifting in $W^{1,p}$  for $p\ge 2$ but not
 for $p<2$:

\begin{theorem} Let $\Omega$ be a bounded simply connected domain in $\mathbb{R}^d,d=2,3$, with  continuous boundary.
 Let $Q\in W^{1,p}(\Omega,\boldsymbol{\mathcal{Q}})$. If $p\ge 2$
there exists a lifting $\varphi^Q\in W^{1,p}(\Omega,\mathbb{S}^2)$ so
that $P\circ\varphi^Q=Q$.
\par Moreover we have the estimate
\begin{equation}
c_1\|\nabla Q\|_{L^p}\le \|\nabla \varphi^Q \|_{L^p}\le c_2\|\nabla
Q\|_{L^p}\label{equiv}
\end{equation} with $c_1,c_2$ constants that depend only on $p$.
\par For $p<2$ there exist line fields for which there is no lifting.  \label{proposition:simplycon}
\end{theorem}
 {\bf Proof.} {\it The case $p\ge 2$.}
\par Let us recall the following  result of  M. Pakhzad and
T. Rivi\`{e}re:

 \medskip\noindent\it {Proposition} {\rm (\cite{PARI03}, p.225)}
Let $M, N$ be compact smooth manifolds, with $M$ simply connected.
For $u\in W^{1,2}(M,N)$ there exists a sequence
$\left\{u^{(k)}\right\}_{k\in\mathbb{N}}$ with $u^{(k)}\in
C^{\infty}(M,N)$ so that $u^{(k)}$ converges weakly to $u$.\rm

\smallskip\par Let us assume first that $\Omega$ is a domain with
smooth boundary. Using the above theorem with $M=\bar\Omega$ and
$N=\boldsymbol{\mathcal{Q}}$
 we can find a sequence of smooth functions $Q^{(k)}$ converging weakly to $Q$.
 Each $Q^{(k)}$ is orientable, by {\rm Proposition}~\ref{continuouslift} as the domain $\Omega$ is simply connected. Using {\rm Proposition}~\ref{limitorient}    we obtain that the
limit function $Q$ is also orientable.
\par We obtain thus that for $Q\in W^{1,p}(\Omega,\boldsymbol{\mathcal{Q}})\subset W^{1,2}(\Omega,\boldsymbol{\mathcal{Q}}),p\ge 2$
there exists $n\in W^{1,2}(\Omega,\mathbb{S}^2)$ with $P(n)=Q$. But
then $n$ is continuous along almost any line parallel with the axis
of coordinates, hence by (\ref{derivtrick}) we get
$n\in W^{1,p}(\Omega,\mathbb{S}^2)$.
\par In order to extend the theorem to less smooth domains we need
the following

\begin{lemma} {\rm (\cite{ballzarnescu+})}
Let $\Omega\subset\mathbb{R}^n$ be a bounded  simply-connected domain with continuous boundary. There exists
$\varepsilon_0>0$ so that for any $\varepsilon>0$ with
$\varepsilon<\varepsilon_0$ there exists a connected open  set
$\Omega_\varepsilon\subset\Omega$ with smooth boundary and such that
$d_H(\Omega_\varepsilon,\Omega)<\varepsilon$ where $d_H$ denotes the Hausdorff distance.
Moreover $\Omega_\varepsilon$ can be chosen so that it is simply
connected and $\Omega_{\varepsilon'}\subset \Omega_\varepsilon$ if $\varepsilon<\varepsilon'$.
\end{lemma}
\par Using the lemma one finds a sequence of simply-connected
smooth domains $\Omega_{\varepsilon_k}\subset
\Omega,k\in\mathbb{N},$ with
$\Omega_{\varepsilon_k}\subset\Omega_{\varepsilon_{k+1}}$ and
$\cup_{k\in\mathbb{N}}\Omega_{\varepsilon_k}=\Omega$. Then for
$\Omega_{\varepsilon_1}$ one has, by the previous arguments, that
there exists $n_{\varepsilon_1}\in
W^{1,p}(\Omega_{\varepsilon_1},\mathbb{S}^2)$ so that
$P(n_{\varepsilon_1})=Q$ on $\Omega_{\varepsilon_1}$. On
$\Omega_{\varepsilon_2}$ one has two possibilities of orienting $Q$,
and one chooses $n_{\varepsilon_2}\in
W^{1,p}(\Omega_{\varepsilon_2},\mathbb{S}^2)$ so that
$n_{\varepsilon_2}(x)=n_{\varepsilon_1}(x)$, a.e.
$x\in\Omega_{\varepsilon_1}$. One continues similarly defining
inductively $n_{\varepsilon_k},k\in\mathbb{N}$.

\par We can define now $n\in W^{1,2}(\Omega,\mathbb{S}^2)$ by $n(x)=n_{\varepsilon_k}(x),\,\mbox{ for all } x\in\Omega_{\varepsilon_k}$.

\par The formula (\ref{equiv}) is  straightforward by taking into
account  the relation between $\varphi^u$ and $u$ as well as
(\ref{derivtrick}).

\smallskip
\noindent {\it The case  $1\le p<2$.}

\par  An example is provided in Figure 2.
\begin{figure}[h]
 \centering
\includegraphics[scale=0.4]{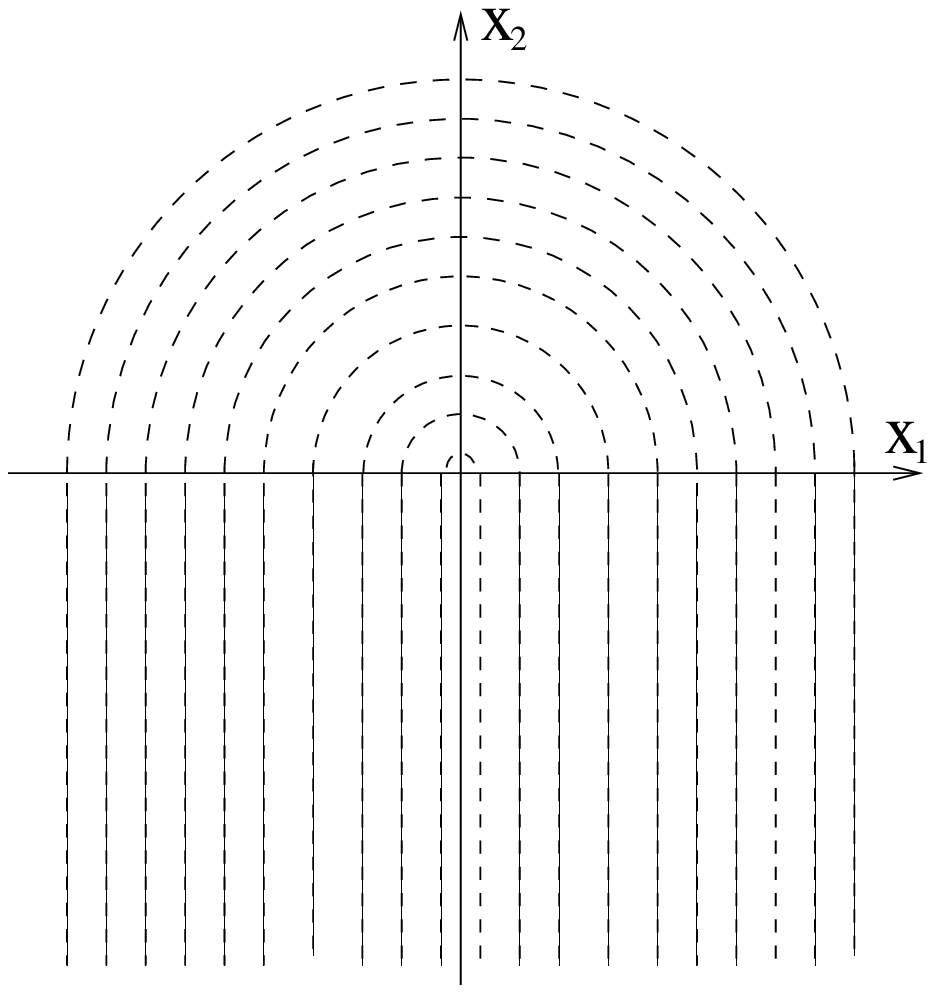}
\caption{A non-orientable director field  on a simply connected
domain, for $p<2$}
\end{figure}
 The line field $$Q(x)=s\left(n(x)\otimes n(x)-\frac{1}{3}Id\right)$$ on $\Omega=(-1,1)^3\subset\mathbb{R}^3$ corresponds to
 what is called in the physical literature `an index one-half singularity', where $$n(x_1,x_2,x_3)=\left\{\begin{array}{ll}
              (x_2,-x_1,0) & \textrm{ for $(x_1,x_2,x_3)\in [0,1)\times (-1,1)^2$}\\
              (0,1,0) &  \textrm{ for  $(x_1,x_2,x_3)\in (-1,0)\times (-1,1)^2$}\\
              \end{array}\right.. $$ $\Box$

\bigskip\section{ Analytic orientability criteria in $2D$}
\label{analytic}
\par In this section we  restrict  ourselves to
 planar line fields, i.e. the domain $\Omega$ is a
subset of $\mathbb{R}^2$ and the line field  takes values only in
$\boldsymbol{\mathcal{Q}}_2$, not in the whole of
$\boldsymbol{\mathcal{Q}}$.

\par  It is important to know,  from a PDE perspective, if it is
possible to detect the orientability (or non-orientability) of a
line field just by knowing its boundary values.

Let us first recall Remark~\ref{boundaryorient:2d} (after
Proposition~\ref{boundaryorient}) which shows
 that orientability in a domain  implies orientability at the
boundary. Thus, in particular, if  a line field in $W^{1,2}(\Omega)$
is orientable then its trace on the boundary,  a line field in
$H^{1/2}(\partial\Omega)$, must be  orientable as well. We will see,
in the next section, in Proposition~\ref{holesorient},
 that the converse is true as well, namely that orientability at the boundary implies orientability in the interior.
  We recall that it was already shown in
Section~\ref{subsection:alternative},
Theorem~\ref{theorem:alternative}, that for continuous line fields
on domains with holes $G$ (as defined in
Section~\ref{subsection:alternative}) orientability can be checked
at the boundary.

\par In order to obtain an analogue of the previous theorem for less regular functions, in $W^{1,2}(G)$, we need first to understand the relation between the orientability in the class of continuous line fields and that in $W^{1,2}(G)$. We study first this question at the boundary and consider a line field on $\partial G$ that is in $C(\partial G)\cap
H^{1/2}(\partial G)$. We claim  that if the line field is
non-orientable, as a continuous line field, then it is also
non-orientable in $H^{1/2}(\partial G)$.

\par More precisely, for $G\subset\mathbb{R}^2$,  a domain with holes as defined in Section~\ref{subsection:alternative},  let $Q\in C(\partial G,\boldsymbol{\mathcal{Q}}_2)\cap
H^{1/2}(\partial G,\boldsymbol{\mathcal{Q}}_2)$ be a line field
non-orientable in the class of continuous fields. We assume for
contradiction that $Q$ is orientable as  a function in
$H^{1/2}(\partial G)$ i.e. that
$Q_{ij}=s(n_in_j-\frac{\delta_{ij}}{3}),i,j=1,2,3$ with $n_i\in
H^{1/2}(\partial G,\mathbb{R}),i=1,2,3$ and
$(n_1(x),n_2(x),n_3(x))\in\mathbb{S}^1$ a.e. $x\in\partial G$. If we
can show that $n_i, i=1,2,3$ has a continuous representative, we
obtain a contradiction which proves our claim.  As $G$ has a smooth
boundary we can, without loss of generality, assume that there
exists locally a smooth transformation that takes functions in
$H^{1/2}(B_\delta(P)\cap \partial G)$ (with $P\in\partial G$) into
functions in $H^{1/2}(I)$ where $I$ is an open interval and thus we
need to show that if $Q_{ii}=sn_in_i\in C(I)$ and $n_i\in
H^{1/2}(I;\mathbb{R})$ for $i=1,2,3$ then there exists a continuous
$\bar n_i,i=1,2,3$ such that $\bar n_i(x)=n_i(x)$ a.e. $x\in I$.

\begin{lemma} Let $I\subset \mathbb{R}$ be an open set. Take $f:I\to \mathbb{R}$ be such that
 $f\in H^{1/2}(I;\mathbb{R})$ and $f^2\in
C(I;\mathbb{R})$. Then there exists $f^*\in C(I,\mathbb{R})$ so that
$f^*=f$ a.e. on $I$. \label{h12continuu}
\end{lemma}
\smallskip{\bf Proof.} We claim first that if $f(a)\ne 0$ there
exists a $\bar\delta=\bar\delta(a)>0$ so that on
$(a-\bar\delta,a+\bar\delta)$ the function $f$ has constant sign
almost everywhere.
\par Assuming the claim the proof is straightforward.  Indeed, let $Z(f)$ denote the zero set of $f$ in $I$. We define $s:I\setminus Z(f)\to \{1,-1\}$ such that $s(y)=1$ if there exists a $\delta_0(y)>0$ so that $f$ is positive almost everywhere on the interval $(y-\delta,y+\delta)$ for any $\delta<\delta_0$, and $s(y)=-1$ otherwise. One can easily check that $s$ is constant on the connected component of $y$  for any $y\in I\setminus Z(f)$.
\par Recalling that $f^2\in C(I;\mathbb{R})$ and so  is defined everywhere we let
$$f^*(y)\stackrel{\rm{def}}{=}\left\{\begin{array}{ll}s(y)\sqrt{f^2} & \textrm{ if }y\in  I\setminus Z(f)\\
                                          0    &\textrm{ if } y\in Z(f)\end{array}\right.$$  and one can easily check that $f^*$ is continuous. Indeed, if $y\in I\setminus Z(f)$ there exists an open interval around $y$, say $(y-\delta,y+\delta)$ on which $s(y)$ is constant hence on $(y-\delta,y+\delta)$ we have that $f^*$ is either plus or minus $\sqrt{f^2}$, and $\sqrt{f^2}$ is a continuous function. If $y\in Z(f)$ let us take $(y_n)_{n\in\mathbb{N}}$  an arbitrary sequence of points so that $y_n\to y$. The continuity of $f^2$ implies that for any $\varepsilon>0$ there exists a $n(\varepsilon)$ such that $|f^2(y_n)|<\varepsilon$ if $n>n(\varepsilon)$ so that $|f(y_n)|=|\pm f^*(y_n)|=|f^*(y_n)|\le \sqrt{\varepsilon}$, which proves the continuity of $f^*$ at $y$.

\smallskip
\par We continue by proving the claim and start by assuming without loss of
generality that $f(a)=l>0$. As $f^2\in C(I,\mathbb{R})$ there exists
a $\delta_0>0$ such that
\begin{equation}
|f^2(x)-l^2|<\frac{l^2}{4},\mbox{ for all } x\in
(a-\delta_0,a+\delta_0). \label{fl}
\end{equation}

\par Note that $H^{1/2}(I,\mathbb{R})\subset
VMO(I,\mathbb{R})$ (see for instance
\cite{BRNI95},\cite{BRNI96},\cite{TRIE95}). Recall that if $f\in
VMO(I)$ then for any $\varepsilon>0$ there exists a $\tilde\delta>0$
such that:
$$\frac{1}{|B(x,\delta)|}\int_{B(x,\delta)}\bigg|f(s)-\frac{1}{|B(x,\delta)|}\int_{B(x,\delta)}f(t)dt\bigg|ds<\varepsilon$$
for all $\delta<\min\{\tilde\delta,\frac{1}{2}(x,\partial I)\}$.

\par We show that there exists a $\delta_1<\tilde\delta$ so that for any
$I(x)_{\delta}=(x-\delta,x+\delta)\subset(a-\delta_0,a+\delta_0)$
with $\delta<\delta_1$ we have

\begin{equation}
\bigg|\frac{1}{|I(x)_\delta|}\int_{I(x)_\delta}
f(y)dy\bigg|>\frac{l}{8}. \label{flpm}
\end{equation}
\par Indeed, if (\ref{flpm}) were false
there would exist two sequences
$(\delta_k)_{k\in\mathbb{N}},(x_k)_{k\in\mathbb{N}}$, with
$\delta_k\to 0$ so that
$I(x_k)_{\delta_k}=(x_k-\delta_k,x_k+\delta_k)\subset(a-\delta_0,a+\delta_0)$
and
 \begin{equation} -\frac{l}{8}\le
\frac{1}{|I(x_k)_{\delta_k}|}\int_{I(x_k)_{\delta_k}}
f(s)ds\le\frac{l}{8}.\label{flpmbis}
\end{equation}

\noindent From the $VMO$ characterization of $f$ we have
$$\frac{1}{|I(x_k)_{\delta_k}|}\int_{I(x_k)_{\delta_k}}\Big|f(s)-\frac{1}{|I(x_k)_{\delta_k}|}\int_{I(x_k)_{\delta_k}}
f(t)dt\Big|ds<\frac{l}{4}$$ for $\delta_k$ small enough. However,
the last inequality cannot hold simultaneously with (\ref{fl}) and
(\ref{flpmbis}). This contradiction proves (\ref{flpm}).

\par Let us denote

$$g(x,\delta)\stackrel{\rm{def}}{=}\frac{1}{|I(x)_\delta|}\int_{I(x)_\delta} f(y)dy.$$

\noindent As $f^2$ (and thus $f$) is bounded on $[a-\delta_0,a+\delta_0]$ one can easily check that $g(x,\delta)$ is continuous as a function of two variables on the set $\{(x,\delta); (x,\delta)\in (a-\delta_0+\delta_1,a+\delta_0-\delta_1)\times [0,\delta_1]\}$ and has no zeros on this set (because of  (\ref{flpm})). Thus $g$ has constant sign on $\{(x,\delta); (x,\delta)\in (a-\delta_0+\delta_1,a+\delta_0-\delta_1)\times [0,\delta_1]\}$ and then by using the Lebesgue differentiation theorem we obtain that  $f$ also has constant sign. $\Box$

\smallskip\par In order to study the orientability of
 planar line fields for $Q$ a $\boldsymbol{\mathcal{Q}}_2$-valued function we
define the {\it auxiliary} complex-valued map
 $A(Q)$:
\begin{equation}
A(Q)\stackrel{\rm{def}}{=}\frac{2}{s}Q_{11}-\frac{1}{3}+i\frac{2}{s}Q_{12},\,\,A(Q)\in
\mathbb{S}^1\subset\mathbb{C}. \label{aux}
\end{equation}

\par The motivation for this definition is that if $Q$ has the form in (\ref{q2})  and $\mathcal{Z}(n)\stackrel{\rm{def}}{=}n_1+in_2$ then $A(Q)=\mathcal{Z}^2(n)$.  The auxiliary map allows one to associate to a planar line field an {\it auxiliary} unit-length vector field. We shall
determine the orientability of the line field in terms of
topological properties of this auxiliary vector field. We provide
first a necessary and sufficient condition for orientability along
the boundary of bounded smooth sets. This does not suffice for a line field to be orientable
 on the whole domain but
provides a necessary condition for it.

\smallskip \par Before stating the orientability criterion, we need to fix some notations about the degree. Let us recall \cite{HIRS76}, pp. $120-130$, that one can define an  integer degree for a  smooth function $f:M\to N$ at a regular value $y=f(x)$ where $M$ and $N$ are  boundaryless, compact and oriented manifolds of the same dimension. We work only with a connected target manifold ($N=\mathbb{S}^1$), so the degree is independent of the regular value chosen \cite{HIRS76}, Lemma $1.4$, p. $124$. In the case when $M$ has several connected components $M_1,\dots, M_k$ we denote $\deg(f,M)=\sum_{i=1}^k \deg(f,M_i)$ where each $M_i$ is given the orientation induced by the inclusion $M_i\hookrightarrow M$. In the case when $M$ is connected and its orientation is the standard one induced from the ambient space we omit the $M$ and simply write $\deg f$.

\par However, sometimes the degree can be defined for functions that are not necessarily smooth. Let us recall Theorem $A.3$
in \cite{BEGP91} that gives a formula for the degree of a
complex-valued function $f\in H^{1/2}(\mathbb{S}^1,\mathbb{S}^1)$,
namely:
\begin{equation}\deg f=\frac{1}{2\pi i}\int_{\mathbb{S}^1}f^{-1}\frac{\partial f}{\partial \theta}\,d\theta.
\label{h12degree}
\end{equation}
(note that the integral is defined in the sense of distributions
since $f^{-1}=\bar f\in H^{\frac{1}{2}}(\mathbb{S}^1,\mathbb{S}^1)$
and $\frac{\partial f}{\partial \theta}\in
H^{-1/2}(\mathbb{S}^1,\mathbb{S}^1)$).

\begin{proposition} Let $\Omega$ be a smooth, bounded domain in
$\mathbb{R}^2$  and let \\ $Q\in
W^{1,2}(\Omega,\boldsymbol{\mathcal{Q}}_2)$. We denote by
$(\partial\Omega)_i,\,i=1,\dots,k$  the connected components of the
boundary.
\par For any $i\in\{1,2,\dots,k\}$ the function ${\rm Tr}\,Q|_{(\partial\Omega)_i}\in
H^{1/2}((\partial\Omega)_i,\boldsymbol{\mathcal{Q}}_2),$ is
orientable (in the space $H^{1/2}$) if and only if
$\deg(A({\rm Tr}\,Q),(\partial\Omega)_i)\in 2\mathbb{Z}$.
Moreover if there exists a unit-length vector field $n$ such that $\mathcal{Z}(n)=n_1+in_2\in
H^{1/2}((\partial\Omega)_i,\mathbb{S}^1)$ and
$P(n)={\rm Tr}\,Q$ a.e. on $(\partial\Omega)_i$ then
$\deg(n,(\partial\Omega)_i)=\frac{1}{2}\deg(A({\rm Tr}\,Q),(\partial\Omega)_i)$.
\label{h12criterion}
\end{proposition}
\smallskip {\bf Proof.} We can regard $\bar\Omega$ as a manifold with boundary
 and then the topological boundary of the set coincides with the boundary
  as a manifold. The boundary  is then again a manifold. More precisely $\partial \Omega$
is a one-dimensional closed manifold without boundary. Taking into
account the classification theorem for one-dimensional manifolds
(see \cite{MILN65}) we have that each connected component of
$\partial\Omega$ is diffeomorphic to $\mathbb{S}^1$. We continue
thus by assuming, without loss of generality, that for
$i\in\{1,2,\dots,k\}$ we have  $(\partial\Omega)_i=\mathbb{S}^1$.
\par It is easily seen that  ${\rm Tr}\,Q\in
H^{1/2}(\mathbb{S}^1,\boldsymbol{\mathcal{Q}_2})$ is orientable if
and only if  for the function $A(\rm{Tr}\,Q)\in
H^{1/2}(\mathbb{S}^1,\mathbb{S}^1)$  there exists a unit-length vector field $n$ such that $\mathcal{Z}(n)\in
H^{1/2}(\mathbb{S}^1,\mathbb{S}^1)$ and $A({\rm Tr}\,Q)=\mathcal{Z}^2(n)$.
\par  We  claim now that a necessary and sufficient condition  for the existence of a unit-length vector field $n$ that  $\mathcal{Z}(n)\in
H^{1/2}(\mathbb{S}^1,\mathbb{S}^1)$ and $A({\rm Tr}\,Q|_{\mathbb{S}^1})=\mathcal{Z}^2(n)$  is $\deg(A({\rm Tr}\,Q),\mathbb{S}^1)\in 2\mathbb{Z}$.
\par We prove first the necessity. It is known (\cite{BEGP91}, p.21) that for any function  $v\in H^{1/2}(\mathbb{S}^1,\mathbb{S}^1)$
 there exists a number $k=\deg\,v\in\mathbb{Z}$ and a unique (up to an integral multiple of $2\pi$)  $V\in
H^{1/2}(\mathbb{S}^1,\mathbb{R})$ so that $v(z)=z^k\cdot e^{iV(z)}$
a.e. $z\in\mathbb{S}^1$. If we assume that
$A({\rm Tr}\,Q|_{\mathbb{S}^1})=\mathcal{Z}^2(n)$ for some  unit-length vector field $n$ with $\mathcal{Z}(n)\in
H^{1/2}(\mathbb{S}^1,\mathbb{S}^1)$ using the quoted result we have
that there exist
$\alpha=\deg\, A({\rm Tr}\,Q|_{\mathbb{S}^1})\in\mathbb{Z}$,
$\beta=\deg\,n\in\mathbb{Z}$ and $g,h\in
H^{1/2}(\mathbb{S}^1,\mathbb{R})$ so that
$A({\rm Tr}\,Q|_{\mathbb{S}^1})(z)=z^\alpha\cdot e^{ig(z)}$ and
$\mathcal{Z}\big(n(z)\big)=z^\beta\cdot e^{ih(z)}$. The equality
$A({\rm Tr}\,Q)|_{\mathbb{S}^1})=\mathcal{Z}^2(n)$ implies that, a.e. on
$\mathbb{S}^1$, one has:
\begin{equation}
z^{\alpha-2\beta}=e^{i(2h-g)}. \label{degrel}
\end{equation}

\par We claim that the last equality implies $\alpha=2\beta$. Indeed,  we have $2h-g\in H^{1/2}(\mathbb{S}^1,\mathbb{R})$ and  thus (see for instance \cite{composition}, Thm. $2$) $e^{i(2h-g)}\in H^{1/2}(\mathbb{S}^1,\mathbb{S}^1)$.
\par Using formula (\ref{h12degree}) we find that  the expression on the  right hand side of (\ref{degrel}) has degree $0$, while the one on the left hand side has degree $\alpha-2\beta\in\mathbb{Z}$, hence our claim.
 In order to prove the sufficiency let us assume that $\deg A({\rm Tr}\,Q)=2k, \,k\in\mathbb{Z}$. Then, by the previously quoted representation formula in (\cite{BEGP91}, p.21) there exists a $W\in H^{1/2}(\mathbb{S}^1,\mathbb{R})$ so that $A({\rm Tr}Q)(z)=z^{2k}e^{iW(z)}$ and thus there exists a vector field $n$ such that $\mathcal{Z}(n)(z)=z^ke^{iW(z)/2}\in H^{1/2}(\mathbb{S}^1,\mathbb{S}^1)$.
\par The same representation formula immediately gives the last part of the Proposition. $\Box$

\smallskip\par We can now provide a necessary and sufficient condition for orientability
{\it on the whole domain}, in the case of a planar domain with
holes.
\begin{proposition}\label{holesorient}
Let $G$ be a planar domain with holes as defined in {\rm (\ref{holes})},
 Section~{\rm \ref{subsection:alternative}}. Assume moreover that $\partial
G$ is smooth.  Let $Q\in W^{1,2}(G,\boldsymbol{\mathcal{Q}}_2)$.
Then $Q$ is orientable if and only if
$$\deg(A({\rm Tr}\,\,Q|_{\partial\Omega}),\partial\Omega)\in
2\mathbb{Z},\,\deg(A({\rm Tr}\,\,Q|_{\partial\omega_i}),\partial\omega_i)\in
2\mathbb{Z},i=1,\dots,n.$$
\end{proposition}
{\bf Proof.} The necessity of the condition is
a consequence of  Proposition ~\ref{h12criterion} together with
 Proposition ~\ref{boundaryorient}. We show the sufficiency. As
$G\subset \mathbb{R}^2$ we have, from \cite{SCUH83}, that there
exists a sequence of  functions $Q_k\in
C^1(\bar G;\boldsymbol{\mathcal{Q}}_2)$ so that $Q_k\to Q$ in
$W^{1,2}(G,\boldsymbol{\mathcal{Q}}_2)$. We show that for $k$ large
enough $Q_k$ is orientable. First let us observe that we have

\begin{lemma} Let $G$ be an open set in $\mathbb{R}^2$. The function $Q\in
W^{1,2}(G;\boldsymbol{\mathcal{Q}}_2)\cap C(\bar
G;\boldsymbol{\mathcal{Q}}_2)$ is orientable as a function in
$W^{1,2}$ if and only if  it is orientable as a continuous function.
\label{equivalencecontinuous}
\end{lemma}
\smallskip{\bf Proof of the lemma.} We assume first that $Q$ is orientable
in $W^{1,2}$ and show that it is orientable in $C$. Let $n\in
W^{1,2}(G;\mathbb{S}^1)$ be such that $P(n)=Q$. Note that this
implies $n_i^2\in W^{1,2}(G;\mathbb{R})\cap C(\bar
G;\mathbb{R}),i=1,2,3$.

  \par  We prove  first that $n_i\in
W^{1,2}(G;\mathbb{R}),i=1,2,3$ and $n_i^2\in C(\bar G;\mathbb{R})$
imply $n_i=n_i^*$ a.e. for some $n_i^*\in C(\bar G;\mathbb{R}),i=1,2,3$. To
prove this  we claim first that:
\begin{eqnarray}
(C)\,\,\textrm{\it If $x_0\in \bar G$ is such  that $n_i^2(x_0)\not=0$ then there is a neighbourhood
of $x_0$} \nonumber\\
\textrm{\it on which $n_i$ has constant sign almost everywhere}\nonumber
\end{eqnarray}
\noindent Assuming $(C)$ it is straightforward to construct $n_i^*$, in a manner nearly identical to the proof of a similar claim in the proof of Lemma~\ref{h12continuu}. We continue by proving the claim $(C)$. Let $l^2\stackrel{\rm def}{=}n_i^2(x_0),l>0$. There
exists $\varepsilon>0$ such that if $|x-x_0|<\varepsilon$ then
$n_i(x)\in (-\frac{5}{4}l,-\frac{3}{4}l)\cup
(\frac{3}{4}l,\frac{5}{4}l)$. From $n_i\in
W^{1,2}(G,\mathbb{R})$ we have that $n_i$ is continuous along almost
all lines parallel with the coordinate axes, in a suitably chosen reference frame. This suffices for
concluding that $n_i$ has constant sign almost everywhere in $\{x\in
G;|x-x_0|<\varepsilon\}$.
\smallskip
\par Assume on the other hand that $Q$ is orientable as a continuous function,
 i.e. there exists a $n\in C(\bar G,\mathbb{S}^1)$ so that $P(n)=v$. Using  Lemma
~\ref{regularitylemma} we have that $n\in W^{1,2}$. $\Box$

 Continuing the proof of the theorem let us recall \cite{BRNI95} that for a unit-length vector field
$n\in H^{1/2}(\mathbb{S}^1,\mathbb{S}^1)$ there exists a $\delta>0$ (depending on $n$) such that for any other unit-length vector field $m\in H^{1/2}(\mathbb{S}^1,\mathbb{S}^1)$ with $\|n-m\|_{BMO}<\delta$ we have that $m$ has the same degree as $n$. Taking into account the relation between the $BMO(\mathbb{S}^1,\mathbb{S}^1)$ and the $H^{1/2}(\mathbb{S}^1,\mathbb{S}^1)$ norms we have that there exists $\delta_0>0$ so that if $\|n-m\|_{H^{1/2}(\mathbb{S}^1,\mathbb{R})}<\delta_0$ then $n$ and $m$ have the same degree.  Thus for $k$ large enough we have that $\deg(A(\textrm{Tr}\,
Q_k|_{\partial\Omega}),\partial\Omega)\in
2\mathbb{Z},\,\deg(A(\textrm{Tr}\,
Q_k|_{\partial\omega_i}),\partial\omega_i)\in
2\mathbb{Z},i=1,\dots,n$.
\par   Proposition~\ref{h12criterion} shows that
$\textrm{Tr}\,Q_k|_{\partial\Omega},\textrm{Tr}\,Q_k|_{\partial\omega_i},i=1,\dots,n$
are orientable in $H^{1/2}$. Using  Lemma~\ref{h12continuu} we have
that
$\textrm{Tr}\,Q_k|_{\partial\Omega},\textrm{Tr}\,Q_k|_{\partial\omega_i},i=1,\dots,n$
are also orientable in the class of continuous functions. Using
Theorem~{theorem:alternative} we have that $Q_k$ is orientable in
the class of continuous functions. Using
Lemma~\ref{equivalencecontinuous} we obtain that for large enough
$k$ the function $Q_k$ is orientable in $W^{1,2}$. Since strong
convergence preserves orientability (see
Proposition~\ref{limitorient}) we conclude that $Q$ is orientable.
$\Box$

\begin{remark}
 It is known (see for instance \cite{ADAMS75}) that for functions with values in $\mathbb{R}^d$ we have $W^{1,2}\setminus C\not=\emptyset$. However one may ask if for functions with values in $\mathbb{S}^1$ the situation is different. This is not the case, as shown by the vector field: $n(x)=(n_1(x),n_2(x),n_3(x))$ with $n_1(x)=\frac{1}{2}\sin\big(\ln\ln(\frac{k}{|x|})\big)$, $n_2(x)=\sqrt{1-n_1(x)^2}$, $n_3(x)=0$, on $D=\{x\in\mathbb{R}^2,|x|\le 1\}$ (we take $k>1$). Then one can easily check that $n\in W^{1,2}(D;\mathbb{S}^1)\setminus C(D;\mathbb{S}^1)$.
\end{remark}

\smallskip\par The previous proposition shows  that we can determine the orientability to computing certain numbers. However, in specific cases, it may be simpler to just use Lemma~\ref{equivalencecontinuous} and check the orientability at the continuous level of regularity, where topological tools can be more efficient.
\par As an example, consider an analytic description of the line field in Figure $1$. Let

\begin{equation}
\tilde Q=s(\tilde n\otimes \tilde n-\frac{1}{3}Id)\in W^{1,2}(\tilde\Omega,\boldsymbol{\mathcal{Q}}_2)
\label{specQ}
\end{equation} where

\begin{eqnarray}
\tilde\Omega\stackrel{\rm{def}}{=}\{(x,y)\in [-1,1]\times [-1,0], \sqrt{x^2+y^2}\ge \frac{1}{2}\}\nonumber\\
\cup \{(x,y); y\ge 0,\frac{1}{2}\le\sqrt{x^2+y^2}\le 1\}
\label{specOmega}
\end{eqnarray} and

\begin{equation}
\tilde n(x,y)=\left\{\begin{array}{ll} (0,1,0) & \textrm{ if }(x,y)\in\left([-1,1]\times [-1,0]\right)\cap\tilde\Omega\\
                                                  (-\frac{y}{\sqrt{x^2+y^2}},\frac{x}{\sqrt{x^2+y^2}},0) & \textrm{ if } y\ge 0,\frac{1}{2}\le\sqrt{x^2+y^2}\le 1
                                          \end{array}\right.
\label{specn}
\end{equation}

\smallskip
\begin{lemma} The line field $\tilde Q$ as  in {\rm (\ref{specQ})}, {\rm (\ref{specn})} on the domain $\tilde\Omega$ as  in
{\rm (\ref{specOmega})} is not orientable in $W^{1,2}(\tilde\Omega;\boldsymbol{\mathcal{Q}}_2)$ or in $C(\tilde\Omega;\boldsymbol{\mathcal{Q}}_2)$.
\label{lemma:specEX}
\end{lemma}
\smallskip{\bf Proof.} Lemma~\ref{equivalencecontinuous} shows that it suffices to prove the non-orientability in the class of continuous line fields. Let us consider the following subsets of $\tilde\Omega$: $\Omega_1\stackrel{\rm{def}}{=}\{(x,y)\in \tilde\Omega, \,y\le 0\}$, $\Omega_2\stackrel{\rm{def}}{=}\{(x,y)\in \tilde\Omega,\, x\le 0\}$, $\Omega_3\stackrel{\rm{def}}{=}\{(x,y)\in \tilde\Omega,\, x\ge 0\}$. We assume for contradiction that the continuous line field is orientable and try to find an orientation. In $\Omega_1$ there are only two possible orientations (see also Proposition~\ref{prop:2orient}), that is all the unit vectors are $(0,1,0)$ or all are $(0,-1,0)$. Let us assume that we pick the orientation $(0,1,0)$. There are two possible orientations in $\Omega_2$ but since $\Omega_1\cap \Omega_2\not=\emptyset$ and we have already chosen an orientation in $\Omega_1$ we can only pick the orientation $(\frac{y}{\sqrt{x^2+y^2}},-\frac{x}{\sqrt{x^2+y^2}})$ in $\Omega_2$. Also there are two possible orientations in $\Omega_3$ but since $\Omega_1\cap \Omega_3\not=\emptyset$ and we have already chosen an orientation in $\Omega_1$ we can only pick the orientation $(-\frac{y}{\sqrt{x^2+y^2}},\frac{x}{\sqrt{x^2+y^2}})$ in $\Omega_3$. Thus on the line $\{(0,y); y\in [\frac{1}{2},1]\}$ we have both the orientation $(-1,0,0)$ and $(1,0,0)$. Similarily, if we start with the other possible orientation in $\Omega_1$ we also reach a contradiction. $\Box$

\section{The minimizing $Q$-harmonic maps versus minimizing harmonic
maps in the plane}

\par We saw in the previous sections that in order to have a geometry in which both orientable and non-orientable energy minimizers exist we need to allow for a domain that is not simply connected. Propositions \ref{h12criterion} and \ref{holesorient} show that if the boundary data on all components of the boundary is orientable then any line field with that boundary data will be orientable. Moreover, if the boundary data on at least one component of the boundary is non-orientable then any line field with that boundary data will necessarily be non-orientable. Thus full knowledge of the boundary data  completely determines the orientability of the line fields with that boundary data. In order to allow for a geometry with both orientable and non-orientable energy minimizers we need to fix orientable boundary data on only one part of the boundary.
\par The simplest situation one could conceive is to consider a domain with one hole. However in such a domain putting orientable boundary data on one component of the boundary  would imply that any line field with that boundary data is orientable (indeed, let $G=\Omega\setminus\overline{\omega}_1$ and $g:\partial\Omega\to
\boldsymbol{\mathcal{Q}}_2$ be orientable, so that degree of $A(g)$ is even; for any $h:\partial\omega_1\to\boldsymbol{\mathcal{Q}}_2$ we need to have $\deg(A(g))+\deg(A(h))=0$, see for instance  \cite{HIRS76}, p. $126$, and hence $h$ is orientable). Thus we need to take at least two holes. If one puts orientable boundary data on two of the components of the boundary, leaving the third component free, a degree argument as before shows that the boundary data on the third component of the boundary must be orientable as well, hence we can only have orientable line fields.

\par Thus we are led to considering the case of a domain with two holes and orientable boundary data on only one connected component of the boundary. Such a situation is presented in Fig.~\ref{fig:4}. More precisely let us consider the domains (for $\delta>1$):

\begin{equation}
\begin{array}{l}
M_1\stackrel{\rm def}{=} \{x=(x_1,x_2)\in\mathbb{R}^2: x_1^2+(x_2-\delta)^2< 1\}\\
M_2\stackrel{\rm def}{=}\{x=(x_1,x_2)\in\mathbb{R}^2: x_1^2+(x_2+\delta)^2< 1\}\\
M_3\stackrel{\rm def}{=}\{x=(x_1,x_2)\in\mathbb{R}^2: |x_1|< 1; |x_2|\le \delta\}\\
M_4\stackrel{\rm def}{=}\{x=(x_1,x_2)\in\mathbb{R}^2: x_1^2+(x_2-\delta)^2\le \frac{1}{2}\}\\
M_5\stackrel{\rm def}{=} \{x=(x_1,x_2)\in\mathbb{R}^2: x_1^2+(x_2+\delta)^2\le\frac{1}{2}\}
\label{M:domains}
\end{array}
\end{equation} and we define the stadium domain:

\begin{equation}M_\delta=M_1\cup M_2\cup M_3\setminus( M_4\cup M_5).\label{Mdelta}
\end{equation}
\par  On the outer boundary we impose as boundary conditions lines tangent to the boundary, which can be oriented clockwise (as shown in $B$ in Fig. ~\ref{fig:4}) or anticlockwise. Thus we have a simple geometry with boundary conditions that allow both orientable and non-orientable line fields.  We compare the minimizers of

$$\mathcal{I}_\delta(Q)=\int_{M_\delta}|\nabla Q(x)|^2\,dx$$ (in $W^{1,2}(M_\delta,\boldsymbol{\mathcal{Q}}_2)$, subject to the indicated line field boundary conditions on the outer boundary) with the minimizers of

$$\mathcal{J}_\delta(n)=2s^2\int_{M_\delta}|\nabla n(x)|^2\,dx$$ (in $W^{1,2}(M_\delta;\mathbb{S}^1)$, subject to tangent vector-field boundary conditions on the outer boundary). Note that $\mathcal{I}_\delta(n)=\mathcal{J}_\delta(Q)$ when $Q$ is orientable.

 We have:

\begin{figure}[t]
\center\includegraphics[scale=0.4]{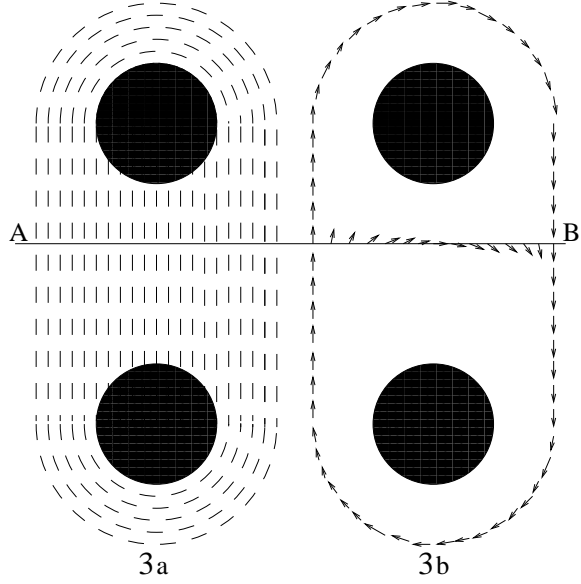} \caption{A situation
in which the energy minimizer is non-orientable} \label{fig:4}
\end{figure}

\begin{lemma}
 Let $\bar n_\delta\in W^{1,2}(M_\delta,\mathbb{S}^1)$ be any global energy minimizer of $\mathcal{I}_\delta(n)$ in $W^{1,2}(M_\delta;\mathbb{S}^1)$ (subject to tangent vector-field boundary conditions on the outer boundary, as in Fig.~\ref{fig:4}b). Let $\bar Q_\delta\in W^{1,2}(M_\delta; \boldsymbol{\mathcal{Q}}_2)$ be any global energy minimizer of $\mathcal{J}_\delta(Q)$ in $W^{1,2}(M_\delta;\boldsymbol{\mathcal{Q}}_2)$ (subject to tangent line-field boundary conditions on the outer boundary, as in Fig.~\ref{fig:4}a).
\par There exists a $\delta_0>1$ so that for any $\delta>\delta_0$ we have
$$\mathcal{J}_\delta(\bar Q_\delta)< \mathcal{I}_\delta(\bar n_\delta).$$

\label{lemma:minimexample}
\end{lemma}
\smallskip{\bf Proof.} Let us observe first that the sets in which we do the minimization, in either the oriented or non-oriented context, are non-empty. Indeed, let us take

$$\tilde Q(x)=\left\{\begin{array}{ll} s\left((0,1,0)\otimes (0,1,0)-\frac{1}{3}Id\right),\, x\in M_3\\
s\left(n_\delta(x)\otimes n_\delta(x)-\frac{1}{3}Id\right),\, x\in M_1\setminus M_4,\, x_2\ge\delta\\
s\left(m_\delta(x)\otimes m_\delta(x)-\frac{1}{3}Id\right),\,x\in M_2\setminus M_5,\, x_2\le -\delta
\end{array}\right.$$ where $$n_\delta(x)\stackrel{\rm def}{=}\left(\frac{x_2-\delta}{|(x_1,x_2-\delta)|},-\frac{x_1}{|(x_1,x_2-\delta)|},0\right)$$
$$m_\delta(x)\stackrel{\rm def}{=}\left(\frac{x_2+\delta}{|(x_1,x_2+\delta)|},-\frac{x_1}{|(x_1,x_2+\delta)|},0\right).$$
\par Then $\tilde Q\in W^{1,2}$ and satisfies the boundary conditions. Let us observe that $\tilde Q$ is exactly the line field shown in Fig.~\ref{fig:4}b. It is also straightforward to see that in the case of  vector-field boundary conditions there exist  vector fields $n_\delta\in W^{1,2}(M_\delta;\mathbb{S}^1)$ on the whole $M_\delta$ that match the boundary conditions.
\par Let us observe that if $n_\delta\in W^{1,2}(M_\delta,\mathbb{S}^1)$ satisfies  the boundary conditions then for almost all $x_2\in [-\delta,\delta]$ we have  $n_\delta(\cdot,x_2)\in W^{1,2}([-1,1];\mathbb{S}^1)$ and $n_\delta(-1,x_2)=(0,1,0)$, $n_\delta(1,x_2)=(0,-1,0)$, and it is an elementary exercise to check that $\int_{[-1,1]\times \{x_2\}} |\partial_{x_1} n_\delta (z,x_2)|^2\,dz\ge \frac{\pi^2}{2}$. Then
\begin{equation}
\int_{M_\delta}|\nabla n_\delta (x)|^2\,dx\ge \int_{M_\delta} |\partial_{x_1} n_\delta (x)|^2\,dx\ge \int_{M_3}|\partial_{x_1}n_\delta (x)|^2\,dx\ge \delta\pi^2.
\end{equation}

\par Thus we have that $\int_{M_\delta}|\nabla \bar n_\delta(x)|^2\,dx\ge \delta\pi^2$ and, noting the way $\tilde Q$ is defined we have that $\int_{M_{\delta}}|\nabla \tilde Q(x)|^2\,dx$ is independent of $\delta$. Hence there exists $\delta_0>0$ so that  for any $\delta>\delta_0$ we have
$$2s^2\int_{M_\delta}|\nabla\bar n_\delta(x)|^2\,dx\ge \delta\pi^2\ge\int_{M_\delta}|\nabla\tilde Q(x)|^2\,dx\ge \int_{M_\delta}|\nabla\bar Q_\delta(x)|^2\,dx$$  which proves the claim. $\Box$
\par The previous theorem shows that for $\delta$ large enough the Oseen-Frank theory  fails to capture the global energy minimizer and  detects just a local energy minimizer, the energy minimizer in the class of oriented line fields. In the following we completely characterize the instances in which the Oseen-Frank theory  fails in this way.

\smallskip
We consider   a smooth planar domain
$G=\Omega\setminus\cup_{i=1}^n\overline{\omega_i}$  with $n\ge 1$
holes, $\omega_i,i=1,\dots,n$, as defined in (\ref{holes}),
Section~\ref{subsection:alternative}. We consider the problem of
minimizing the energy
\begin{equation}\label{zz}
\mathcal{I}_G(Q)=\int_G |\nabla Q(x)|^2\,dx,
\end{equation}
on this domain in the class of
$\boldsymbol{\mathcal{Q}}_2$-valued functions whose gradients are
square integrable and that satisfy $Q|_{\partial\Omega}=g$ with $g$ smooth.  We  shall provide necessary and
sufficient conditions for the global minimizers to be
non-orientable. This is the most interesting situation as it is
precisely that in which the Oseen-Frank theory would fail
to see the right energy minimizer and would only provide a local
energy minimizer, a minimizer in the class of orientable line
fields.

\begin{figure}[h]
 \centering
\includegraphics[scale=0.4]{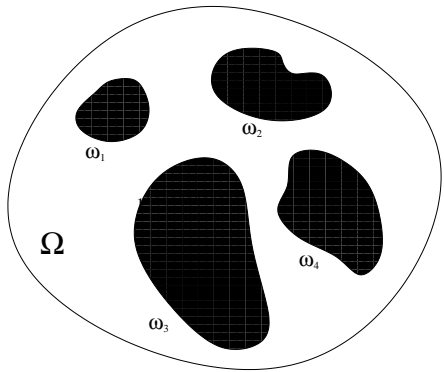}
\caption{A domain with holes}
\end{figure}

\par In order to encode the complexity of the domain and its
relationship with the prescribed boundary data $g$, we need $n+1$
functions $h_1,\dots, h_n$ and $h(g)$.The  functions  $h_i,i=1,2,\dots,n$ encode the characteristics of the holes and their relations with the set $\Omega$. Each function $h_i,i=1,\dots, n$  is the solution of
the equation

\begin{equation}
\left\{\begin{array}{ll} \Delta h_i=0 & \textrm{ on G }\\
                         h_i=1        & \textrm{ on
                         $\partial\omega_i$}\\
                         h_i=0        & \textrm{ on
                         $\partial\omega_j,\,j\not=i$}\\
                         \frac{\partial h_i}{\partial\nu}=0 &\textrm{ on $\partial\Omega$}
\end{array}\right.
\label{hi}
\end{equation}
\par    We define the matrix $D=(D_{ij}),i,j=1,\dots,n$ {\it depending only on the domain}, by $D_{ij}\stackrel{\rm def}{=}\frac{1}{2\pi}\int_{\partial\omega_i} \frac{\partial
h_j}{\partial\nu}(\sigma)\,d\sigma$. Note that $D_{ij}=\frac{1}{2\pi}\int_G \nabla h_i(x)\cdot\nabla h_j(x)\,dx$ so that $D$ is symmetric.
\par It will be important, in  later calculations, to   know explicitly the nullspace of the matrix $D$:

\begin{lemma} Let $e\stackrel{\rm def}{=}(\underbrace{1,1,\dots,1}_{\textrm{ n times}})$. Then the nullspace of $D$ is $N(D)=\mathbb{R}e$.
\label{Dkernel}
\end{lemma}
{\bf Proof.} Let $z\in\mathbb{R}^n$, $h(x)\stackrel{\rm def}{=}(h_1(x),\dots,h_n(x))$ and denote $v(x)\stackrel{\rm def}{=}h(x)\cdot z$. Then
\begin{displaymath}
Dz=0\Leftrightarrow \int_{\partial\omega_i}\frac{\partial v}{\partial \nu}\,d\sigma=0,\,\,\mbox{ for all } i\in\{1,2,\dots,n\}.
\end{displaymath}
\par The last relation implies that $\int_G |\nabla v(x)|^2\,dx=\int_{\partial G} v\cdot\frac{\partial v}{\partial \nu}\,d\sigma=0$ and hence $v$ is a constant function. But $v=z_i$ on $\partial\omega_i$ and so $z_1=z_2=\dots=z_n$.
\par Conversely if $z_1=\dots=z_n=a$ then
\begin{displaymath}
\left\{\begin{array}{ll} \Delta v=0 & \textrm{ on }G\\
 v=a &\textrm{ on }\partial\omega_i\\
 \frac{\partial v}{\partial n}=0 &\textrm{ on each }\partial\Omega
 \end{array}\right.
\end{displaymath} and so by uniqueness $v\equiv a$ and $Dz=0$. $\Box$

\par In order to define the function $h(g)$ we need to use the auxiliary  vector-field $A(g)=\frac{2}{s}g_{11}-\frac{1}{3}+i\frac{2}{s}g_{12}$ associated to the line field $g$,
as defined in (\ref{aux}). This is a complex-valued function but from now on, until the end of the paper, we identify in a standard way the complex-valued function $A(g)$ with a vector-valued real function. The function $h(g)$ describes the relation between the domain and
the boundary data and  is defined as the solution of the equation

\begin{equation}
\left\{\begin{array}{ll} \Delta h(g)=0 & \textrm{ on G }\\
                         \frac{\partial h(g)}{\partial \nu}=A(g)\times\frac{\partial  A(g)}{\partial \tau}
                          & \textrm{ on $\partial\Omega$}\\
                         h(g)=0        & \textrm{ on
                         $\partial G\setminus\partial\Omega$}\\
\end{array}\right.
\label{h0}
\end{equation} (note that in the above the vector-valued real function $A(g)\times\frac{\partial  A(g)}{\partial \tau}$ is identified, in a standard way, with a scalar, real-valued function). The derivative $\frac{\partial}{\partial
\tau}$ is the tangential derivative on the boundary.

\par We  define the vector $J(g)\stackrel{\rm{def}}{=}(J(g)^1,\dots, J(g)^n)$ {(depending on both the domain and  the boundary data)}, where $J(g)^i=\frac{1}{2\pi}\int_{\partial\omega_i} \frac{\partial h(g)}{\partial\nu}\,d\sigma$.

\par Let
$$\mathcal{D}(g)\stackrel{\rm def}{=}\{(d_1,\dots,d_n)\in\mathbb{Z}^n, \sum_{i=1}^n
d_i=-\deg(A(g),\partial\Omega)\}$$
$$\mathcal{D}_{even}(g)\stackrel{\rm def}{=}\{(d_1,\dots,d_n)\in \left(2\mathbb{Z}\right)^n, \sum_{i=1}^n
d_i=-\deg(A(g),\partial\Omega)\}.$$

\bigskip

\par We can now state  a necessary and sufficient criterion for determining the
orientability of the global minimizer of the $Q$-harmonic maps
problem:

\begin{theorem}  Let $g\in W^{1,2}(\partial\Omega,\boldsymbol{\mathcal{Q}}_2)$  be  orientable, and assume that
$$W_g^{1,2}(G,\boldsymbol{\mathcal{Q}}_2)=\{Q:G\to \boldsymbol{\mathcal{Q}}_2;\,\mathcal{I}_G(Q)<\infty,\,Q|_{\partial \Omega}=g\}$$
is nonempty. Then the infimum of $\mathcal{I}_G(Q)$ in $W_g^{1,2}(G,\boldsymbol{\mathcal{Q}}_2)$ is attained.
\par  For $d\in\mathcal{D}(g)$ let  $c(d)\stackrel{\rm def}{=}(c_1(d),\dots,c_n(d))$ be a solution of the equation
\begin{equation}
D\cdot c=d-J(g).
\label{dcequation}
\end{equation}
 Then a necessary and sufficient condition for all global minimizers to be non-orientable is
\begin{equation}
 \inf_{d\in \mathcal{D}(g)} c(d)\cdot Dc(d)<\inf_{d\in \mathcal{D}_{even}(g)}c(d)\cdot Dc(d).\label{condition}
\end{equation}
\label{criterionprop}
\end{theorem}
\smallskip{\bf Proof.} For any $Q\in W_g^{1,2}(G,\boldsymbol{\mathcal{Q}}_2)$
 we have $A(Q)\in \mathbb{S}^1$
and moreover
$$\|\nabla Q\|_{L^2(G)}=\frac{\sqrt{2}}{s}\|\nabla A(Q)\|_{L^2(G)}.$$
 Observing that $A$ is a bijective operator we have that our
minimization problem reduces to

\begin{equation}
\inf_{m\in W_{A(g)}^{1,2}(G,\mathbb{S}^1)} \frac{2}{s^2}\int_G
|\nabla m(x)|^2\,dx. \label{minequiv}
\end{equation}
It is well known that the minimum of the energy for the last
 problem is attained by a function $m_{min}\in W_{A(g)}^{1,2}(G,\mathbb{S}^1)$
 satisfying a harmonic map equation (see \cite{BEBH94}).
 \par In order to determine this function we  first claim that
  for any  $m\in W_{A(g)}^{1,2}(G,\mathbb{S}^1)$ we have
\begin{equation}
\deg(A(g),\partial\Omega)=-\Sigma_{i=1}^n \deg(m,\partial\omega_i).
\label{degsum}
\end{equation}
Indeed, by \cite{SCUH83}  there exists a sequence
 $m_k\in W_{A(g)}^{1,2}(G,\mathbb{S}^1)\cap C(G,\mathbb{S}^1)$ so
 that $m_k\to m$ in $W^{1,2}$. Taking into account that the function $m_k$ is continuous, and the
 properties of the degree for continuous functions, \cite{HIRS76}, p. $126$
   we have
\begin{displaymath}
\deg(A(g),\partial\Omega)=-\sum_{i=1}^n \deg(m_k,\partial\omega_i).
\end{displaymath}
 Using the continuity of the trace operator we
let $k\to\infty$ in the last relation and we obtain the claimed
relation (\ref{degsum}). \par  Thus we can divide the function
space $W_{A(g)}^{1,2}$ into countably many disjoint subsets
corresponding to maps with given degrees $d_i$ on each
$\partial\omega_i,i=1,\dots,n$. A way of solving the minimization
problem (\ref{minequiv}) is to obtain first the minimizer  on each
such subset as before and thus obtain countably many functions
$m_1,m_2,\dots, m_l,\dots,l\in\mathbb{N}$. The solution of
(\ref{minequiv}) is then that $m_k$ with $\|\nabla
m_k\|_{L^2(G)}=\inf_{i\in\mathbb{N}}\|\nabla m_i\|_{L^2(G)}$ (such
an $m_k$ exists because there exists a global minimizer for the
problem (\ref{minequiv})).
\par Thus we need to study first the minimization problem
\begin{equation}
\inf_{\stackrel{m\in
W_{A(g)}^{1,2}(G,\mathbb{S}^1)}{\deg(m,\partial\omega_i)=d_i,i=1,\dots,n}}
\frac{2}{s^2}\int_G |\nabla m(x)|^2\,dx\label{minequiv+}
\end{equation} for each set of $d_i,i=1,\dots,n$ so that
$d=(d_1,\dots,d_n)\in\mathcal{D}(g)$.

\par The advantage in studying (\ref{minequiv+}) rather than (\ref{minequiv}) is that
the determination of the minimum for (\ref{minequiv+}) can be
reduced to a simpler, scalar, problem.   Indeed, it is shown in
\cite{BEBH94} that for a fixed $d\in\mathcal{D}(g)$ if we denote by
$\rm m^*$ a minimizer of (\ref{minequiv+}) then $\|\nabla
{\rm m^*}\|_{L^2(G)}=\|\nabla \Phi\|_{L^2(G)}$ where $\Phi$ is the
unique solution (up to an additive constant) of the scalar problem:
 \begin{equation}
 \left\{\begin{array}{ll}
 \Delta\Phi=0 &\textrm{ on G}\\
 \Phi=c_i     &\textrm{ on $\partial\omega_i$}\\
 \int_{\partial \omega_i}\frac{\partial \Phi}{\partial\nu}\,d\sigma=2\pi d_i
 &\\
 \frac{\partial\Phi}{\partial\nu}=A(g)\times \frac{\partial
 A(g)}{\partial\tau} & \textrm{ on $\partial\Omega$}\end{array}\right.
 \label{eqreduced}
 \end{equation} where $d_i$ are prescribed (with $d=(d_1,\dots, d_n)\in\mathcal{D}(g)$), but not the $c_i$.
 \par Taking $\Phi$  to be a solution of the above problem let us denote
 $h=\Phi-\sum_{i=1}^n c_ih_i$. We have that $h$ is a solution of the
 problem (\ref{h0}). Since (\ref{h0}) has a unique solution we get
 that $h\equiv h(g)$, thus
 \begin{equation}
 \Phi=h(g)+\sum_{i=1}^n c_i h_i.
 \label{phirep}
 \end{equation}
 Taking into account the
 representation (\ref{phirep}) of $\Phi$ as well as the equation it satisfies,
  (\ref{eqreduced}), we get
  \begin{equation}
  \sum_{j=1}^n  D_{ij}c_j+J^i(g)=d_i,i=1,\dots
  n.\label{csystem}
  \end{equation}
  Lemma~\ref{Dkernel} shows that the last system has a one-dimensional affine space of solutions, so that the function $c(d)$ introduced in the statement is multivalued. We claim that, however, the value of $c(d)\cdot Dc(d)$ is independent of the particular representative of $c(d)$ used.
Indeed observe that multiplying  (\ref{dcequation}) by $e$ we obtain $Dc\cdot e=d\cdot e-J(g)\cdot e$ and since $D$ is symmetric and $De=0$ we obtain $d\cdot e=J(g)\cdot e$, thus proving our claim.
 \smallskip

  \par In order to finish the proof it suffices to  recall the orientability criterion given by
   Proposition~\ref{holesorient} and  observe that from
  (\ref{phirep}) we have 
  \begin{eqnarray}
 \|\nabla\Phi\|_{L^2(G)}^2&=&\sum_{i,j=1}^n c_ic_j\int_G \nabla
  h_i\cdot\nabla h_j\,dx+2\sum_{j=1}^n c_j\int_G \nabla h_j\cdot\nabla h(g)\,dx
\nonumber\\&& \hspace{2.5in}+\int_G|\nabla h(g)|^2\,dx\nonumber\\
&\hspace{-.8in}=&\hspace{-.4in}\sum_{i,j=1}^n c_i c_j \int_{\partial G}\frac{\partial h_i}{\partial \nu}h_j\,d\sigma+2\sum_{j=1}^n c_j\int_{\partial G} \frac{\partial h_j}{\partial \nu} h(g)\,d\sigma
+\int_G|\nabla h(g)|^2\,dx\nonumber\\
&\hspace{-.8in}=&\hspace{-.4in}c(d)\cdot Dc(d) +\int_G|\nabla h(g)|^2\,dx,\nonumber
\end{eqnarray} where we used the definitions of $h(g)$,\,$h_i,i=1,\dots,n$ and $D_{ij},i,j=1,\dots,n$.
 $\Box$
\smallskip\begin{remark} One can see, carefully following the proof, that one also has that if
\begin{equation}
 \inf_{d\in \mathcal{D}(g)\setminus  \mathcal{D}_{even}(g)} c(d)\cdot Dc(d)>\inf_{d\in \mathcal{D}_{even}(g)}c(d)\cdot Dc(d)\label{condition+}
\end{equation} then all global energy minimizers must necessarily be orientable.
\par Moreover, if
\begin{equation}
\inf_{d\in \mathcal{D}(g)\setminus  \mathcal{D}_{even}(g)} c(d)\cdot Dc(d)=\inf_{d\in \mathcal{D}_{even}(g)}c(d)\cdot Dc(d)\label{equalitycondition}
\end{equation} then there exist both an orientable and a non-orientable global energy minimizer.
\label{remark:orientablemin}
\end{remark}

\smallskip\begin{remark} One can easily see that for $g$ smooth the set $W^{1,2}_g$ is non-empty by recalling \cite{guillemin}, \cite{HIRS76} that a degree zero smooth map $\tilde g:\partial G\to\mathbb{S}^1$ can be extended to a smooth map on $G$. On the other hand one can alway choose some suitable smooth vector field $h:\partial G\setminus \partial\Omega\to\mathbb{S}^1$ so that $$\tilde g(x)=\left\{\begin{array}{ll} g(x) &\textrm{ if } x\in \partial\Omega\\
h(x) &\textrm{ if }x\in\partial G\setminus\partial\Omega\end{array} \right.$$ has degree zero.
\par In general, for $g$ not smooth, the space $W^{1,2}_g$ may be empty (see  \cite{partialconstrainedbdry}).
\label{remark:nonemptyfunctionalspace}
\end{remark}

\bigskip We continue with a more in-depth analysis of the case when the domain $G$ has only two holes, by using the tools developed in the previous Proposition.

\begin{proposition}{\rm(i)} Let $G=\Omega\setminus\cup_{i=1}^2\overline{\omega_i}$ with two holes, $\omega_1$ and $\omega_2$, as defined in {\rm (\ref{holes})}, Section~{\rm \ref{subsection:alternative}}. We take a boundary data $g\in W^{1,2}(\partial\Omega,\boldsymbol{\mathcal{Q}}_2)$ that is an orientable line field and assume that the space $W^{1,2}_g(\Omega,\boldsymbol{\mathcal{Q}}_2)$ is non-empty. Then $\rm{dist} (J(g)^1,\mathbb{Z})=\rm{dist} (J(g)^2,\mathbb{Z})$, $\rm{dist} (J(g)^1,2\mathbb{Z})=\rm{dist} (J(g)^2,2\mathbb{Z})
$ and  all the global energy minimizers are non-orientable if and only if

\begin{equation}
\rm{dist} (J(g)^1,\mathbb{Z})<\rm{dist} (J(g)^1,2\mathbb{Z}).\label{2holescriterion}
\end{equation}
\par On the other hand, if

\begin{equation}
\rm{dist}(J(g)^1,2\mathbb{Z})<\rm{dist}(J(g)^1,2\mathbb{Z}+1)
\label{criterion:orientablemin}
\end{equation} then all the global energy minimizers are orientable.
\par Moreover, if
\begin{equation}
\rm{dist}(J(g)^1,2\mathbb{Z})=\rm{dist}(J(g)^1,2\mathbb{Z}+1)
\label{criterion:orientablemin2}
\end{equation} then there exist both an orientable and a non-orientable energy minimizer.

\smallskip\noindent
 {\rm (ii)} Let $M_\delta$ be the domain defined in {\rm(\ref{Mdelta})}.  Let $\bar n_\delta\in W^{1,2}(M_\delta,\mathbb{S}^1)$ be any global energy minimizer of $\mathcal{I}_\delta(n)$ in $W^{1,2}(M_\delta;\mathbb{S}^1)$ (subject to tangent vector-field boundary conditions on the outer boundary, as in Fig.~{\rm \ref{fig:4}b)}. Let $\bar Q_\delta\in W^{1,2}(M_\delta; \boldsymbol{\mathcal{Q}}_2)$ be any global energy minimizer of $\mathcal{J}_\delta(Q)$ in $W^{1,2}(M_\delta;\boldsymbol{\mathcal{Q}}_2)$ $($subject to tangent line-field boundary conditions on the outer boundary, as in Fig.~{\rm\ref{fig:4}a}$)$.
\par For any $\delta>1$ we have
$$\mathcal{J}_\delta(\bar Q_\delta)< \mathcal{I}_\delta(\bar n_\delta).$$

\smallskip\noindent
{\rm (iii)} Let $$G_\delta\stackrel{\rm def}{=}\{x=(x_1,x_2)\in\mathbb{R}^2: \frac{1}{2}<x_1^2+x_2^2<1,\, x_1^2+(x_2-\frac{3}{4})^2>\delta\}$$ for $\delta<\frac{1}{4}$. Let $g:\{(x_1,x_2)\in\mathbb{R}^2:x_1^2+x_2^2=1\}\to \boldsymbol{\mathcal{Q}}_2$ be a smooth orientable line field and let $\tilde g$ to be a vector field that is one of the orientations of $g$.
\par There exists a $\delta_0$ so that for any $\delta<\delta_0$ any global energy minimizer of $I_\delta(Q)$ must necessarily be orientable.

\smallskip\noindent
{\rm (iv)}  Let
$$M_6^\rho\stackrel{\rm def}{=}\{x=(x_1,x_2)\in\mathbb{R}^2: x_1^2+(x_2-2)^2\le \rho\}$$ for $\rho<1$. With $M_i,\,i=1,\dots,5$ defined as in {\rm (\ref{M:domains})} {\rm(}where we take $\delta=2${\rm)} we consider the domain:
$$N_\rho\stackrel{\rm def}{=}M_1\cup M_2\cup M_3\setminus (M_5\cup M_6^\rho)$$
\par We impose tangential line field boundary conditions on the outer boundary of $N_\rho$. Then there exists a $\rho\in (0,\frac{1}{2})$ so that there exist both an orientable and a non-orientable  global energy minimizer of $I_{N_\rho}(Q)$, subject to the imposed boundary conditions.
\label{proposition:2holes}
\end{proposition}
{\bf Proof.} {\rm (i)}  Let $d_1,d_2\in\mathbb{Z}$ be some arbitrary pair such that
\begin{equation}
d_1+d_2=-\deg (A(g),\partial\Omega).
\label{d1d2}
\end{equation}

\par Relation $d\cdot e=J(g)\cdot e$ (with $e=(1,1)$), shown in the proof of Theorem~\ref{criterionprop}, gives that:
\begin{equation}
d_1+d_2=J(g)^1+J(g)^2.
\label{djgequality}
\end{equation}
\par Corresponding to this domain $G$ we have the symmetric matrix $D=\left(\begin{array}{ll} a & b\\ b & c\end{array}\right)$ and $D(1,1)^t=0$. Thus $b=-a,\,a=c$ and $$D=a\left(\begin{array}{ll} \,\,1 &  -1\\  -1 & \,\,\,1 \end{array}\right).$$
\par So  equation (\ref{dcequation}) becomes
\begin{equation}
Dc=a\left(\begin{array}{l} c_1-c_2\\  c_2-c_1\end{array}\right)=\left(\begin{array}{l} d_1\\  d_2 \end{array}\right)-\left(\begin{array}{l} J(g)^1\\  J(g)^2 \end{array}\right).
\label{dcequation+}
\end{equation}
\par We claim now that $a\not=0$. Assuming for contradiction that $a=0$ equation (\ref{dcequation+}) (which always has a solution by the arguments in the proof of Proposition~\ref{criterionprop})  implies $d_1=J(g)^1$ and $d_2=J(g)^2$. However if we   replace $d_1$ by $d_1+1$ and $d_2$ by $d_2-1$ (so that their sum is still $-\deg(A(g),\partial\Omega)$) then  (\ref{dcequation+})  no longer has a solution. Thus our claim is proved.
\par Moreover, recalling that $D_{ij}=\frac{1}{2\pi}\int_G \nabla h_i(x)\nabla h_j(x)\,dx$ one can easily see that $D$ is non-negative definite, so $a>0$. Hence \begin{eqnarray}
c\cdot Dc&=&(d_1-J(g)^1)c_1+(d_2-J(g)^2)c_2\nonumber\\
&\stackrel{\rm (\ref{djgequality})}{=}&(d_1-J(g)^1)(c_1-c_2)\stackrel{\rm (\ref{dcequation+})}{=}\frac{1}{a}(d_1-J(g)^1)^2.\nonumber\end{eqnarray}
\par Thus Proposition~\ref{criterionprop} implies that all minimizers are non-orientable if and only if (\ref{2holescriterion}) holds.  Relations (\ref{d1d2}) and (\ref{djgequality}) together with the assumption that $g$ is orientable, hence $\deg(A(g),\partial\Omega)$ is even means that   we always have $\rm{dist}(J(g)^1,2\mathbb{Z})=\rm{dist}(J(g)^2,2\mathbb{Z})$ and $\rm{dist}(J(g)^1,2\mathbb{Z}+1)=\rm{dist}(J(g)^2,2\mathbb{Z}+1)$.

\par The claimed criteria are now a consequence of Remark~\ref{remark:orientablemin}.

\smallskip\noindent {\rm (ii)} We claim that the symmetry of the domain and that of the boundary data imply that $h(g)(x_1,x_2)=h(g)(x_1,-x_2)$. Indeed, let $\tilde h(g)(x_1,x_2)\stackrel{\rm def}{=} h(g)(x_1,-x_2).$

\par  One can check that $$A(g)\times \frac{\partial A(g)}{\partial\tau}(x_1,x_2)=\left\{\begin{array}{ll} 2 &\textrm{ if } (x_1,x_2)\in\partial(\cup_{i=1}^3 M_i),\, |x_2|\ge \delta\\
0 &\textrm{ if }(x_1,x_2)\in\partial(\cup_{i=1}^3 M_i),\,|x_2|<\delta\end{array}\right.$$  and thus the boundary data is symmetric with respect to the $x_2=0$ axis.
Then $h(g)$ and $\tilde h(g)$ are both  functions that solve problem (\ref{h0}), that has a unique solution. Thus $h(g)=\tilde h(g)$ and our claim is proved.
\par Let us observe that the line-field example in the proof of Lemma~\ref{lemma:minimexample} shows that $W_g^{1,2}\not=\emptyset$. Taking this into account, together with  $g\in W^{1,2}(\partial\Omega)$,  we can use the first part of the Proposition (in our case the domain in only $C^1$ and not smooth but one can check that $C^1$ regularity suffices for using the Theorem~\ref{criterionprop} and thus the first part of the Proposition). The symmetry of $h(g)$ implies that $J(g)^1=J(g)^2$ and relations  (\ref{d1d2}), (\ref{djgequality}) together with $\deg(A(g),\partial (M_1\cup M_2\cup M_3))= 2$ imply $J(g)^1=J(g)^2=-1$. Hence the criterion (\ref{2holescriterion}) holds for any $\delta>1$.

\smallskip\noindent {\rm (iii)} Let us first observe that since we took the boundary data to be smooth the function space $W^{1,2}_g(G_\delta,\boldsymbol{\mathcal{Q}}_2)$ is non-empty (see Remark~\ref{remark:nonemptyfunctionalspace}). Let $\Xi_\delta\stackrel{\rm def}{=}\{x=(x_1,x_2)\in\mathbb{R}^2: x_1^2+(x_2-\frac{3}{4})^2<\delta^2\}$ which we regard as the hole $\omega_1$ and  denote by $h^\delta(g)$ the solution of (\ref{h0}) for the domain $G_\delta$. We claim that
\begin{equation}
\int_{\partial\Xi_{\frac{1}{8}}} \frac{\partial h^\delta(g)}{\partial \nu}\,d\sigma\to 0\,\textrm{ as }\delta\to 0.
\label{limithdelta}
\end{equation}

\par On the other hand the divergence theorem shows that  (for $\delta<\frac{1}{8}$):

\begin{displaymath}
0=\int_{\Xi_{\frac{1}{8}}\setminus\Xi_\delta} \Delta h^\delta(g)\,dx=\int_{\partial\Xi_{\frac{1}{8}}}\frac{\partial h^\delta(g)}{\partial\nu}\,d\sigma-\underbrace{\int_{\partial\Xi_\delta}\frac{\partial h^\delta(g)}{\partial\nu}\,d\sigma}_{=2\pi\cdot J(g)^1}.
\end{displaymath}

\par The last relation, together with our previous claim and the criterion (\ref{criterion:orientablemin})  finish the proof.

\smallskip
\par It remains to prove our claim (\ref{limithdelta}) and to this end let us consider the solution $H$ of the pure Neumann problem:
\begin{equation}
\left\{\begin{array}{ll}
\Delta H=0 &\textrm{ on }B_1(0)\setminus B_{1/2}(0)\\
\frac{\partial H}{\partial \nu}=A(g)\times \frac{\partial A(g)}{\partial\tau} &\textrm{ on }\partial B_1(0)\\
\frac{\partial H}{\partial\nu}(\frac{1}{2},\theta)=-2\Big(A(g)\times \frac{\partial A(g)}{\partial\tau}\Big)(1,\theta) & \textrm{ for }\theta\in [0,2\pi]\\
\int_{B_1(0)\setminus B_{\frac{1}{2}}(0)} H(x)\,dx=0
\end{array}\right.
\label{pureN}
\end{equation}
 where on the third line of the system above we used polar coordinates $(r,\theta)\in [\frac{1}{2},1]\times [0,2\pi]$.

\par The solution is smooth on $B_1(0)\setminus B_{1/2}(0)$ and continuous on the closure, \cite{lionsmagenes}, thus there exist  $c_1,c_2>0$ so that $c_2>H(x)>-c_1$ for all $x\in\overline{B_1(0)\setminus B_{1/2}(0)}$. Let $w(x)\stackrel{\rm def}{=}h^\delta(g)(x)-c_1-H(x)$. Then $\Delta w=0$ on the set $G_\delta$. As $\frac{\partial w}{\partial\nu}=0$ on $\partial B_1(0)$ Hopf's lemma shows that $w$ cannot attain its maximum on $\partial B_1(0)$.  Hence, by the maximum principle, it attains its maximum on $\partial G_\delta\setminus\partial B_1(0)$ where, by our construction, $w\le 0$. Thus $w\le 0$ on $\overline{G_\delta}$ and hence $h^\delta(g)\le c_1+H$ on $\overline{G_\delta}\subset\overline{B_1(0)\setminus B_{1/2}(0)}$. Similarily, taking the function $v=H-c_2-h^\delta(g)$ and reasoning analogously we obtain that $H-c_2\le h^\delta(g)$ on $\overline{G_\delta}\subset\overline{ B_1(0)\setminus B_{1/2}(0)}$. Thus

\begin{equation}
H-c_2\le h^\delta(g)\le c_1+H, \,\,\textrm{ on }\overline{G_\delta}\subset\overline{ B_1(0)\setminus B_{1/2}(0)}
\label{uniformhg}
\end{equation}
 and since the sequence of harmonic functions $h^\delta (g)$ is uniformly bounded on a sequence of domains shrinking into the annulus $B_1(0)\setminus B_{1/2}(0)$, we obtain \cite{gilbargtrudinger}, p.$23$, by taking a diagonal sequence, that there exists a function  $f$  so that $h^{\delta_j}(g)$ converges uniformly on compact subsets of  $\Big(B_1(0)\setminus B_{1/2}(0)\Big)\setminus \{(0,\frac{3}{4})\}$ to $f$. Then $f$ is harmonic  \cite{gilbargtrudinger}, p.$23$ on   $\Big(B_1(0)\setminus B_{1/2}(0)\Big)\setminus \{(0,\frac{3}{4})\}$ and bounded (by (\ref{uniformhg})) and hence it has a removable singularity at $(0,\frac{3}{4})$ and  is harmonic on $B_1(0)\setminus B_{1/2}(0)$. Then $\int_{\partial\Xi_{\frac{1}{8}}}\frac{\partial f}{\partial \nu}\,d\sigma=\int_{\Xi_{\frac{1}{8}}}\Delta f(x)\,dx=0$ and since $h(g)^{(\delta)}$ converges uniformly (and thus in $C^\infty$ since $h(g)^{(\delta)}$ are harmonic) on compact sets to $f$ we obtain the claimed relation (\ref{limithdelta}).

\smallskip\noindent {\rm (iv)} Let $g$ correspond to tangential boundary conditions on the outer boundary of $N_\rho$. The example constructed in Lemma~\ref{lemma:minimexample} can be easily modified to show that $W^{1,2}_g(N_\rho)\not=\emptyset$, for all $\rho\in (0,1)$. We denote by $S\stackrel{\rm def}{=}\cup_{i=1}^3 M_i$ the stadium without holes.

 \par Let  $\tilde H_0^1(N_\rho)\stackrel{\rm def}{=}\{ u\in H^1(N_\rho); \textrm{Tr}\,u\Big|_{\partial M_5\cup\partial M_6^\rho}=0\}$. Let $\varphi\in C^\infty(N_\rho)\cap H_0^1(N_\rho)$  be a function vanishing in a neighbourhood of $\partial M_5\cup\partial M_6^\rho$ and denote by $\tilde \varphi$ its extension by zero to a function on $S$.  Denoting by $ h(g)^\rho\in \tilde H^1_0(N_\rho)$ the solution of problem (\ref{h0}) on $N_\rho$,  we  have:
 \begin{equation}
 \int_{N_{\rho}}\nabla h(g)^\rho\cdot\nabla\varphi\,dx=\int_{\partial N_\rho} \frac{\partial h(g)^\rho}{\partial \nu}\varphi\,d\sigma=\int_{\partial S}A(g)\times \frac{\partial A(g)}{\partial\tau}\varphi\,d\sigma
 \label{weakh0}
 \end{equation} and then
 \begin{eqnarray}&\hspace{-.5in}
\left|\int_{N_\rho}\nabla h(g)^\rho\cdot\nabla\varphi\,dx\right| =\left|\int_{\partial S}\frac{\partial h(g)^\rho}{\partial \nu}\varphi\,d\sigma\right|\nonumber\\
& \le \|\frac{\partial h(g)^\rho}{\partial \nu}\|_{H^{-1/2}(\partial S)}\|\varphi\|_{H^{1/2}(\partial S)} \le C_1\|\frac{\partial h(g)^\rho}{\partial \nu}\|_{H^{-1/2}(\partial S)}\|\tilde \varphi\|_{H^1(S)}.
 \nonumber
 \end{eqnarray}
 In the  last inequality  we can assume without loss of generality that $C_1$ is a constant independent of $\rho$, because  $\varphi=\tilde\varphi $ on  $\partial S$ and the last inequality expresses the continuity of the trace operator in $H^1(S)$.
\par We denote by $\tilde h(g)^\rho$ the extension by zero of $h(g)^\rho$ to a function on $S$  and then the last inequality implies
\begin{equation}
\left|\int_{S}\nabla\tilde h(g)^\rho\cdot\nabla\tilde \varphi\,dx\right|\le C_1\|\frac{\partial h(g)^\rho}{\partial \nu}\|_{H^{-1/2}(\partial S)}\|\tilde \varphi\|_{H^1(S)}.
\end{equation}
 Replacing $\tilde\varphi$ in the inequality by $\tilde\varphi_k$  with  $\tilde \varphi_k\to\tilde h(g)^\rho$ in $H^1(S)$ as $k\to\infty$  the last inequality implies:
\begin{equation}
\|\nabla\tilde h(g)^\rho\|_{L^2(S)}^2\le C_1\|A(g)\times \frac{\partial
 A(g)}{\partial\tau}\|_{H^{-1/2}(\partial S)}\|\tilde h(g)^\rho\|_{H^1(S)}.
\label{halfh1unifobds}
\end{equation}
\par We  denote $J(g)^2_\rho\stackrel{\rm def}{=}\frac{1}{2\pi}\int_{\partial M_5} \frac{\partial h(g)^\rho}{\partial \nu}\,ds$ and claim that
\begin{equation}
\rho\mapsto J(g)^2_\rho,\,(\rho\in (0,\frac{1}{2}]),\textrm{ is a continuous function}.
\label{continuityclaim}
\end{equation}
 In order to prove the claim we argue by contradiction and assume that there exists a $\bar\rho\in (0,\frac{1}{2}],\varepsilon_0>0$ and a sequence $\rho_k$ with $\rho_k\to\bar \rho$ and
\begin{equation}
|J(g)^2_{\rho_k}-J(g)^2_{\bar \rho}|>\varepsilon_0,\mbox{ for all } k.
\label{contrassumpt}
\end{equation}
 First let us observe that for $\rho\in (\frac{\bar\rho}{2},\frac{1}{2})$ the functions $\tilde h(g)^\rho$ are zero on a common set of non-zero measure $M_6^{\frac{\bar\rho}{2}}$. Then one has a Poincar\'e inequality:
\begin{equation}
\|\tilde h(g)^\rho\|_{L^2(S)}\le C_2\|\nabla\tilde h(g)^\rho\|_{L^2(S)}
\label{l2unifobds}
\end{equation} for $\rho\in (\frac{\bar\rho}{2},\frac{1}{2})$, with $C_2$  depending on $S$ and $\bar\rho$ (see the argument in \cite{ziemer}, p. $177$ that can be checked to hold even for $p=n=2$).

\par Relations (\ref{halfh1unifobds}) and (\ref{l2unifobds}) imply that there exists a subsequence of $\rho_k$, relabelled as the initial sequence, such that

\begin{equation}
\tilde h(g)^{\rho_k}\rightharpoonup L\,\textrm{ in }H^1(S)
\label{weakconvhg}
\end{equation} for some function $L\in H^1(S)$. We claim that $L\equiv\tilde h(g)^{\bar \rho}$. To this end it suffices to show that for any smooth $\hat\varphi\in H^1(N_{\bar\rho})$  vanishing in a neighbourhood of  $ M_6^{\bar\rho}\cup M_5$ we have:
\begin{equation}
\int_{S}\nabla L\cdot\nabla\hat\varphi\,dx=\int_{\partial S} A(g)\times \frac{\partial A(g)}{\partial \tau}\hat\varphi\,d\sigma
\label{weakhgbarrho}
\end{equation} which, by the uniqueness of the weak solution for the problem (\ref{h0}), implies the claim.

\par In order to prove the last equality let us take a sequence $\varphi^k\in H^1(S)\cap C^\infty(S)$, s
upported in $N_{\rho_k}$ and vanishing near $\partial
M_6^{\rho_k}\cup\partial M_5$, so that
$\varphi^k\to\hat\varphi$ in
$H^1(S)$. Then replacing $N_\rho$ by $S$, $h(g)^\rho$ by $\tilde
h(g)^{\rho_k}$, $\varphi$ by $\varphi^k$ in (\ref{weakh0}) and
passing to the limit $k\to\infty$ by using (\ref{weakconvhg}), we
obtain relation (\ref{weakhgbarrho}).

\smallskip

\par Using that $\nabla h(g)^\rho\in H(W,\textrm{div})$ with $W\subset S\setminus (M_5\cup M_6^{2\bar\rho})$
an open set such that $\partial M_5\subset \partial W$,  relation (\ref{weakconvhg}) and the continuity of the
 normal part of the trace in the space $H(W,\textrm{div})$ (see \cite{lionsmagenes}) we have that  $\frac{\partial h(g)^{\rho_k}}{\partial\nu}\rightharpoonup\frac{\partial h(g)^{\bar\rho}}{\partial\nu}$ in $H^{-\frac{1}{2}}(\partial M_5)$. On the other hand we can write $J(g)^2_{\rho_k}=\frac{1}{2\pi}\int_{\partial M_5} \frac{\partial h(g)^{\rho_k}}{\partial\nu}\,d\sigma=\frac{1}{2\pi}\int_{\partial M_5}\langle \frac{\partial h(g)^{\rho_k}}{\partial\nu},1\rangle\,d\sigma$, with $\langle\cdot,\cdot\rangle$ denoting the duality between $H^{-\frac{1}{2}}(\partial M_5)$ and $H^{\frac{1}{2}}(\partial M_5)$,  and the previously proved weak convergence implies $J(g)^2_{\rho_k}\to J(g)^2_{\bar\rho}$ thus contradicting (\ref{contrassumpt}). The contradiction we have reached proves our earlier claim (\ref{continuityclaim}).
\par The argument in part $(iii)$ can be easily adapted to show that $J(g)^2_\rho\to 0$ as $\rho\to 0$ and the
 argument in part $(ii)$ shows that $J(g)^2_{\frac{1}{2}}=-1$. Thus (\ref{continuityclaim}) shows there exists
  some $\rho_0\in (0,\frac{1}{2})$ so that $J(g)^2_{\rho_0}=-\frac{1}{2}$. Then the criterion (\ref{criterion:orientablemin2})
   shows that on $N_{\rho_0}$ there exist both an orientable and a non-orientable energy minimizer. $\Box$

 \smallskip\begin{remark} Part {\rm (ii)} of  Proposition~\ref{proposition:2holes} shows that in Lemma~\ref{lemma:minimexample}, one does not in fact need the assumption that $\delta$ is large enough.
 \end{remark}
\smallskip\begin{remark} The proof of Part {\rm (iv)} of  Proposition~\ref{proposition:2holes} can be easily adapted to show that for any domain with two holes, if one shrinks enough one of the holes (while keeping the other hole and the orientable boundary data unchanged) then the global energy minimizer will necessarily be orientable.
\end{remark}
\rm{\bf Acknowledgement}
 JMB is very grateful to Apala Majumdar for introducing him to the Landau~-~de Gennes theory; some of the questions addressed here had their origins in the proposal for the EPSRC grant EP/E010288/1 which was written jointly with her and which supported the research  of both authors together with the EPSRC Science and Innovation award to the Oxford Centre for Nonlinear PDE (EP/E035027/1). AZ thanks Luc Nguyen and Yves Capdeboscq for useful discussions concerning regularity issues. We also thank Ulrike Tillmann for useful comments.

\bigskip

\appendix
\section{The relation between sets of class $C^k$ (Lipschitz) and $C^k$ (Lipschitz) manifolds with boundary}

\par We recall (see \cite{MILN65}, p.12, \cite{HIRS76}, pp. 29-30) that a subset $X\subset \mathbb{R}^m$  is  an $m$-dimensional  $C^k$(Lipschitz) manifold with boundary, for $k\ge 0$, if each $x\in X$ has a neighbourhood $U$ in $\mathbb{R}^m$ such that $U\cap X$  is $C^k$-diffeomorphic (homeomorphic for $k=0$, or bi-Lipschitz homeomorphic in the case of Lipschitz manifolds) with a set $V\cap H^m$ in $H^m$, where $V$ is an open set in $\mathbb{R}^m$ and  $$H^m=\{(x_1,\dots,x_m)\in\mathbb{R}^m\,|x_m\ge 0\}.$$
\par  The {\it boundary} of the manifold $X$ is the set of all points in $X$ which correspond to points of $\partial H\stackrel{\rm{def}}{=}\mathbb{R}^{m-1}\times\{0\}\subset \mathbb{R}^m$ under such a diffeomorphism. It can be shown that the condition for a point in $X$ to be on the boundary is independent of the chart chosen (see \cite{HIRS76}).
\par In order to show that a domain $\Omega$ of class $C^k$ (Lipschitz) is a manifold with boundary and the topological boundary coincides with the boundary as a manifold it suffices to show that, for $P\in\partial\Omega$ and $\delta>0$ such that $B_\delta(P)\cap \partial\Omega$ can be represented as the graph of a function, there exists a $C^k$-diffeomorphism $F$ from $B_\delta(P)\cap\overline{\Omega}$ to $V\cap H^m$ which carries $B_\delta(P)\cap\partial\Omega$ into $\partial H^m\cap V$.
\par To this end we let $f:\mathbb{R}^{m-1}\times \{0\}\subset \mathbb{R}^m\to \mathbb{R}$  be such that

\begin{equation}
B_\delta(P)\cap\Omega=\{y=(y',y_m)\in\mathbb{R}^m:\,y_m>f(y',0),\,|y'|<\delta\}\nonumber
\end{equation} where $f$ is of class $C^k$. Then it can be easily checked that  $F:B_{\frac{\delta}{2}}(P)\cap\overline{\Omega}\to V\cap H^m$ defined by
$F(y',y_m)=(y',y_m-f(y',0))$ is injective, onto  the image (we take $V$ so that $F$ is onto) and a $C^k$-diffeomorphism (respectively bi-Lipschitz)  with inverse $F^{-1}(y',y_m)=(y',y_m+f(y',0))$. Moreover $F$ carries $B_\delta (P)\cap\partial\Omega$ bijectively into $V\cap\partial H_m$.

\bibliographystyle{plain}
\bibliography{Orientability}

\end{document}